\documentclass[preprint,12pt]{elsarticle}

\usepackage{amssymb}

\usepackage[activate={true,nocompatibility},final,tracking=true,kerning=true,spacing=true,factor=1100,stretch=10,shrink=10]{microtype}
\usepackage[utf8]{inputenc}
\usepackage[T1]{fontenc}
\usepackage{babel}

\usepackage{graphicx,color}
\usepackage{amsmath}
\usepackage[version=4]{mhchem}
\usepackage{siunitx}
\usepackage{longtable,tabularx}
\usepackage{subcaption}
\usepackage{hyperref}
\usepackage{ragged2e}
\usepackage{arydshln}
\usepackage{optidef} 
\usepackage[ruled,vlined]{algorithm2e}  
\usepackage{mathrsfs}   
\usepackage{changepage}   
\usepackage{float} 
\usepackage{dsfont}
\usepackage{nccmath}
\usepackage{tabularx}
\usepackage{booktabs}
\usepackage{longtable}
\usepackage{enumitem}
\usepackage{adjustbox}
\usepackage{amsmath}               
\usepackage{amsfonts}               
\usepackage{stanli}
\newcommand{\eqnref}[1]{Eq.~(\ref{#1})}
\newcommand{\figref}[1]{Fig.~\ref{#1}}

\newcommand{\tabref}[1]{Table~\ref{#1}}

\usepackage{bm}

\usepackage{optidef} 
\usepackage[ruled,vlined]{algorithm2e}  
\usepackage{mathrsfs}   
\usepackage{changepage}   
\usepackage{float} 
\usepackage{dsfont}

\usepackage{amsthm}
\usepackage{relsize}

\usepackage{nomencl}
\usepackage{multicol}
\renewcommand{\nompreamble}{\begin{center}\begin{multicols}{2}}
\renewcommand{\nompostamble}{\end{multicols}\end{center}}
\makenomenclature
\providetoggle{nomsort}
\settoggle{nomsort}{true} 
\makeatletter
\iftoggle{nomsort}{%
    \let\old@@@nomenclature=\@@@nomenclature        
        \newcounter{@nomcount} \setcounter{@nomcount}{0}%
        \renewcommand\the@nomcount{\two@digits{\value{@nomcount}}}
        \def\@@@nomenclature[#1]#2#3{
          \addtocounter{@nomcount}{1}%
        \def\@tempa{#2}\def\@tempb{#3}%
          \protected@write\@nomenclaturefile{}%
          {\string\nomenclatureentry{\the@nomcount\nom@verb\@tempa @[{\nom@verb\@tempa}]%
          \begingroup\nom@verb\@tempb\protect\nomeqref{\theequation}%
          |nompageref}{\thepage}}%
          \endgroup
          \@esphack}%
      }{}
\makeatother      

\def\abs#1{\left\lvert#1\right\rvert}

\usepackage{xpatch}
\makeatletter
\xpatchcmd\SetKwInOut
  {\hangafter=1\parbox[t]}
  {\hangafter=1\justify\parbox[t]}
  {}{\fail}
\makeatother

\SetKwInOut{Input}{Input}
\DontPrintSemicolon

\usepackage{comment}

\newcommand{\tstar}[5]{
\pgfmathsetmacro{\starangle}{360/#3}
\draw[#5] (#4:#1)
\foreach \x in {1,...,#3}
{ -- (#4+\x*\starangle-\starangle/2:#2) -- (#4+\x*\starangle:#1)
}
-- cycle;
}

\journal{Neurocomputing. DOI 10.1016/j.neucom.2023.126472}

\begin{document}

\begin{frontmatter}

\newtheorem{assumption}{Assumption}
\newtheorem{theorem}{Theorem}
\newtheorem{corollary}{Corollary}
\newtheorem{lemma}{Lemma}
\newtheorem{remark}{Remark}

\title{A mixed-categorical correlation kernel for Gaussian process}
\author[First,Second]{P. Saves\corref{cor1}}
\ead{paul.saves@onera.fr}
\author[Third]{Y. Diouane}
\ead{youssef.diouane@polymtl.ca}
\author[First]{N. Bartoli}
\ead{nathalie.bartoli@onera.fr}
\author[First]{T. Lefebvre}
\ead{thierry.lefebvre@onera.fr}
\author[Fourth]{J. Morlier}
\ead{joseph.morlier@isae-supaero.fr}

\cortext[cor1]{Corresponding author.}

\affiliation[First]{organization={ONERA/DTIS},
            addressline={Universit{\'e} de Toulouse}, 
            city={F-31055 Toulouse},
            country={France}}

\affiliation[Second]{organization={ISAE-SUPAERO},
              addressline={Universit{\'e} de Toulouse}, 
            city={Toulouse},
            postcode={31055 Cedex 4}, 
            country={France}}

\affiliation[Third]{ organization={Polytechnique Montr{\'e}al},
            city={Montr{\'e}al},
            state = {Qu{\'e}bec},
            country={Canada}}

\affiliation[Fourth]{ organization={Institut Cl{\'e}ment Ader (ICA), Universit{\'e} de Toulouse,   ISAE-SUPAERO,  Mines Albi, UPS, INSA, CNRS},
             adressline={3 rue Caroline Aigle},
            city={31400 Toulouse},
            country={France}}

\begin{abstract}
Recently, there has been a growing interest for mixed-categorical meta-models based on Gaussian process (GP) surrogates. In this setting, several existing approaches use different strategies either by using continuous kernels (\textit{e.g.}, continuous relaxation and Gower distance based GP) or by using a direct estimation of the correlation matrix. 
In this paper, we present a kernel-based approach that extends continuous exponential kernels to handle mixed-categorical variables. The proposed kernel leads to a new GP surrogate that generalizes both the continuous relaxation and the Gower distance based GP models. 
We demonstrate, on both analytical and engineering problems, that our proposed GP model gives a higher likelihood and a smaller residual error than the other kernel-based state-of-the-art models. Our method is available in the open-source software SMT. 
\end{abstract}

\end{frontmatter}

\section{Introduction}
\renewcommand*\footnoterule{}

\label{sec:intro}

Expensive-to-evaluate blackbox simulations play a key role for many engineering and industrial applications. In this context, surrogate models have shown great interest for a wide range of applications, \textit{e.g.}, aircraft design~\cite{SciTech_cat}, deep neural networks~\cite{snoek2015scalable}, coastal flooding prediction~\cite{lopez}, agriculture forecasting~\cite{MLP}, turtle retinas modeling~\cite{retina} or seismic imaging~\cite{YDiouane_SGratton_XVasseur_LNVicente_HCalandra_2016}.
These blackbox simulations are generally complex and may involve mixed-categorical input variables. Typically, an aircraft design tool has to take into account variables such as the number of panels, the list of cross sectional areas or the material choices.

In this work, we target to learn an inexpensive surrogate model $ \hat{f}$ from a mixed-categorical blackbox function given by
\begin{equation}
f :  \Omega \times S \times \mathbb{F}^l \to \mathbb{R}.
  \label{eq:opt_prob}
\end{equation}
This function $f$ is typically an expensive-to-evaluate simulation with no exploitable derivative information.
 $\Omega \subset \mathbb{R}^n$ represents the bounded continuous design set for the $n$ continuous variables.  $S \subset \mathbb{Z}^m$ represents the bounded integer set where $L_1, \ldots, L_m$ are the numbers of levels of the $m$ quantitative integer variables on which we can define an order relation and $ \mathbb{F}^l = \{1, \ldots, L_1\} \times \{1, \ldots, L_2\} \times  \ldots \times \{1, \ldots, L_l\}$ is the design space for the $l$ categorical qualitative variables with their respective  $L_1, \ldots, L_l$ levels.
\textcolor{black}{Typical examples of $f$ can be found in different  engineering contexts. Mechanical performance of hybrid discontinuous composite materials~\cite{RaulAIAA} is an example where the mixed-categorical function $f$ represents the stiffness value which depends on a set of input variables $w=(x,c)\in \Omega \times \mathbb{F}^2$. The continuous part $x$ has two components, the length of the fibers $x_1$ and the proportion of carbon fibers $x_2$ (i.e., $\Omega =
[515,12000] \times [0,1]$). The categorical choices $c$ represent the types of carbon fibers $c_1$ and glass ones  $c_2$ (i.e., $\mathbb{F}^2 = \{\text{XN-90},\text{T800H}\} \times \{  \text{GF}, \text{T300},\text{C100},\text{C320} \} $).} 

For that purpose, \textit{Gaussian process} (GP)~\cite{williams2006gaussian}, also called Kriging model~\cite{krige1951statistical}, is known to be a good modeling strategy to learn a response surface model from a given dataset. Namely,
we will consider that our unknown blackbox function $f$ follows a Gaussian process of mean $\mu^{f}$ and of standard deviation $\sigma^f$, $\textit{i.e.}$,
\begin{equation}
f \sim \hat{f}=\mbox{GP}\left(\mu^{f}, [\sigma^f]^2\right). \label{eq:GP:f}\end{equation}
For a general problem involving categorical or integer variables, several modeling strategies to build a mixed-categorical GP have been proposed~\cite{Pelamatti, Zhou, Deng, Roustant,GMHL,Gower,cuesta2021comparison,SciTech_cat}. Compared to a continuous GP, the major changes are in the estimation of the correlation matrix, the latter being essential to build estimates of $\mu^{f}$ and $\sigma^f$. Similarly to the process of constructing a GP with continuous inputs, relaxation techniques~\cite{GMHL,SciTech_cat}, continuous latent variables~\cite{cuesta2021comparison} and Gower distance based models~\cite{Gower} use a kernel-based approach to estimate the correlation matrix.
Other recent approaches try to estimate the correlation matrix independently of a kernel choice by modeling  directly the possible correlation entries of the correlation matrix~\cite{Pelamatti, Zhou, Deng, Roustant}. 

Using GP surrogates is not the only possible approach whatsoever. Random forests are often used instead of GP as they also can model both mean and variance~\cite{SMAC} and tree-structured Parzen estimators have been shown to be well-adapted for such problems~\cite{TPE}. Other surrogate models for blackbox include ReLU functions~\cite{Relu-surr}, piecewise linear neural network~\cite{nn-surr} or categorical regression splines~\cite{splines}. Models other than GP could also be based on a mixed integer kernel as for support vector regression~\cite{herrera} or on a mixed integer distance as for radial basis functions~\cite{RBF_geo}. Another classical modeling strategy is to consider a different continuous model for every possible categorical choice and to build another model peculiar to the categorical variables besides the continuous models. This categorical model can be, for instance, a probability law~\cite{CAT-EGO}, a multi-arm bandit~\cite{Bandit-BO} or an integer model~\cite{AMIEGO}. Also, in case of prior information, latent variables approaches~\cite{cuesta2021comparison} and user-defined neighbourhood~\cite{Mixed_Abramson} based models are of great interest.

In this paper, we target to extend the classical paradigm for continuous inputs (where a kernel is used to build the GP) to cover the mixed-categorical case. Namely, we will present a kernel-based approach that will lead to a unified model for existing approximation strategies~\cite{Pelamatti, GMHL,Gower}. Namely, this work unifies both distance based kernels and matrix based kernels into a unique homogeneous formulation.
This work generalizes existing methods that were already proven to be efficient over deep learning models~\cite{GMHL} and analytical test cases~\cite{Gower}. 
A similar kernel for the estimation of the correlation matrix could be applied to continuous, integer and categorical inputs. The good potential of the proposed approach is shown and analyzed over analytical and industrial test cases. 

\textcolor{black}{Another main benefit behind the use of specific kernels~\cite{Pelamatti, Zhou, Deng, Roustant} to handle mixed-categorical inputs is to model accurately correlations between the  variables; which is required to get accurate GP models. It might be possible to use continuous kernels to model categorical data but in this case one needs to define a distance function. Such function is not trivial to define on categorical data; only simple distances are possible (e.g., Gower distance~\cite{Gower}) which in general leads to poor GP models. This paper shows in particular the utility of mixed-categorical kernel over continuous based ones on both analytical and industrial test cases.
}

The GP models and the Bayesian Optimization (BO) that could be performed with them are implemented in the Surrogate Modeling Toolbox (SMT) v2.0\footnote{\url{https://smt.readthedocs.io/en/latest/}}~\cite{SMT2019}.
Our modeling software is free and open-source and has been used regularly in the aircraft industry, for example with a deep learning model~\cite{DL1, DL2, DL3, DL4} or with a deep gaussian process~\cite{DGP1,DGP2}. 

The remainder of this paper is as follows.
In Section~\ref{sec:GP}, a detailed review of the GP model for continuous and for categorical inputs is given.
The extended kernel-based approach for constructing the correlation matrix is presented in Section~\ref{sec:model_uni}.
Section~\ref{sec:Results} presents academical tests as well as the obtained results.
Conclusions and perspectives are finally drawn in Section~\ref{sec:conclu}.

\section{GP for mixed-categorical inputs}
\renewcommand{\footnoterule}{ \hrule width5cm \vspace*{0.1cm} }    
\label{sec:GP}

In this section, we will present the mathematical background associated with GP for mixed-categorical variables. This part also introduces the notations that will be used throughout the paper. In this section, we are considering the general case involving mixed integer variables. Namely, we assume that $f:\mathbb{R}^n \times  \mathbb{Z}^m \times \mathbb{F}^l \mapsto \mathbb{R}$ and our goal is to build a GP surrogate model for $f$. 

Given a set of data points, called a Design of Experiments (DoE)~\cite{forrester}, Bayesian inference learns the GP model that explains the best this data set. A GP model consists of a mean response hypersurface $\mu^{f}$, as well as an estimation of its variance $[\sigma^f]^2$. In the following, $n_t$ denotes the size of the given DoE data set $(W, \textbf{y}^f)$ such that $W=\{w^1,w^2,\ldots,w^{n_t}\} \in (\mathbb{R}^n \times  \mathbb{Z}^m \times \mathbb{F}^l)^{n_t}$ and $\textbf{y}^f=[f(w^1),f(w^2),\ldots,f(w^{n_t})]^{\top}$. 
For an arbitrary $w= (x,z,c) \in \mathbb{R}^n \times  \mathbb{Z}^m \times \mathbb{F}^l$, not necessary in the DoE, the GP model prediction at $w$ writes as $\hat f(w) = \mu ({w})+\epsilon(w) \in \mathbb{R} $, with $\epsilon$ being the error between $f$ and the model approximation $\mu$~\cite{GP14}. The considered error terms are random variables of variance $\sigma^2$.  Using the DoE, the expression of  $\mu^{f}$ and the estimation of its variance $[\sigma^f]^2$ are given as follows:

\begin{equation} \label{eq:mean:GP}
\mu^f(w)= \hat{\mu}^f+r(w)^\top  [R(\Theta)]^{-1}(\textbf{y}^f-\mathds{1} \hat{\mu}^f), 
\end{equation}
and
\begin{equation}
\label{eq:std:GP}
[\sigma^f(w)]^2=[\hat{\sigma}^f]^2\left[1-r(w)^\top  [R(\Theta)]^{-1}r(w)+ \frac{ \left(1-\mathds{1}^\top  [R(\Theta)]^{-1}r(w) \right)^2}{\mathds{1}^\top  [R(\Theta)]^{-1}\mathds{1}}\right], \end{equation}
where $\hat{\mu}^f$ and $\hat{\sigma}^f$, respectively, are the maximum likelihood estimator (MLE)~\cite{MLE} of $\mu$ and $\sigma$. $\mathds{1}$ denotes the vector of $n_t$ ones. $R$ is the $ n_t \times n_t $ correlation matrix between the input points and $r(w)$ is the correlation vector between the input points and a given $w$.
The correlation matrix $R$ is defined, for a given couple  $(r,s) \in (\{1,\ldots,n_t\})^2$, by \begin{equation}
\label{eq:R}
 [R(\Theta)]_{r,s}=k\left(w^r,w^s,\Theta\right) \in \mathbb{R},\end{equation}
and the vector $r(w)\in \mathbb{R}^{n_t}$ is defined as $ r(w) =[k(w,w^1), \ldots , k(w,w^{n_t})]^{\top}$,
where $k$ is a given correlation kernel that relies on \textcolor{black}{a set of hyperparameters $\Theta$~\cite{de2015optimizing,bernal2022criteria}}. 
The mixed-categorical  correlation kernel is given as the product of three kernels:
\begin{equation}
k(w^r,w^s,\Theta) =  k^{cont}\left(x^r,x^s,\theta^{cont}\right) k^{int}\left(z^r,z^s,\theta^{int}\right)
k^{cat}\left(c^r,c^s,\theta^{cat}\right),
\label{eq:decomp_mix}
\end{equation}
where $k^{cont}$  and $\theta^{cont}$ are the continuous kernel and its associated hyperparameters, $k^{int}$  and $\theta^{int}$ are the integer kernel and its hyperparameters, and last $k^{cat}$  and $\theta^{cat}$ are the ones related with the categorical inputs. In this case, one has $\Theta=\{ \theta^{cont},\theta^{int},\theta^{cat}\}$. 
Henceforth, the general correlation matrix $R$ will rely only on the set of the hyperparameters $\Theta$:
\begin{equation}
    \label{eq:corel:mat}
    [R(\Theta)]_{r,s} = [R^{cont}(\theta^{cont})]_{r,s}  [R^{int}(\theta^{int})]_{r,s}
    [R^{cat}(\theta^{cat})]_{r,s},
\end{equation}
where $[R^{cont}(\theta^{cont})]_{r,s} =k^{cont}(x^r,x^s,\theta^{cont}) $, $[R^{int}(\theta^{int})]_{r,s} =k^{int}(z^r,z^s,\theta^{int}) $ and  $[R^{cat}(\theta^{cat})]_{r,s}=k^{cat}(c^r,c^s,\theta^{cat})$. 
The set of hyperparameters $\Theta$ could be estimated using the DoE data set $({W},\textbf{y}^f)$ through the MLE approach on the following way
\begin{equation}
\Theta^*= \arg\max_{\Theta} \mathcal{L}(\Theta):=\left( - \frac{1}{2} {\textbf{y}^f}^\top [R(\Theta)]^{-1} {\textbf{y}^f}   - \frac{1}{2} \log 	\abs{  [R(\Theta)]} - \frac{n_t}{2} \log 2 \pi    \right),
\label{eq:likelihood}
\end{equation}
where $R(\Theta)$ is computed using~\eqnref{eq:corel:mat}. 
To construct the correlation matrix, several choices for the correlation kernel are possible. Usual families of kernels include exponential kernels or Matern kernels~\cite{Lee2011}. 
In the rest of this section, we will focus mainly on the exponential kernels and describe in details the construction of the continuous $R^{cont}(\theta^{cont})$, the integer $R^{int}(\theta^{int})$ and the categorical $R^{cat}(\theta^{cat})$ correlation matrices. 

\subsection{Correlation matrices for continuous and integer inputs}
\label{sec:GP_cont}

The construction of the correlation matrix $R^{cont}(\theta^{cont})$ for continuous inputs, based on an exponential kernel, can be described as follows. 
For a couple of continuous inputs $x^r\in \mathbb{R}^n$ and $x^s\in \mathbb{R}^n$, one sets 
\begin{equation}
[R^{cont}(\theta^{cont})]_{r,s}= {\displaystyle \prod_{j=1}^{n}{  \exp \left(-   \theta^{cont}_j  \abs{x^r_j - x^s_j}^p \right) }}. 
\label{eq:E_cont}
\end{equation}
Different values for $p$ can be used. Typically, when $p=1$, one gets the absolute exponential kernel (Ornstein-Uhlenbeck process~\cite{Lee2011}) and, when $p=2$, the squared exponential kernel (or Gaussian kernel~\cite{williams2006gaussian}) is obtained. Clearly, in the continuous case, constructing $R^{cont}(\theta^{cont})$ would require the estimation of $n$ non-negative hyperparameters, $\textit{i.e.}$, $\theta^{cont} \in \mathbb{R}^n_+$.

Thanks to a continuous relaxation technique that transforms integer inputs into continuous ones, the integer inputs can be naturally handled with continuous kernels. On this base, in what comes next, there will be no distinction between continuous and integer inputs; the two of them will be handled in the same way. In fact, for integer variables, the distance defined in the continuous case is still valid. Thus, for an integer couple $z^r\in \mathbb{Z}^m$ and $z^s\in \mathbb{Z}^m$, a natural extension of the exponential kernel that handles integer variables can be given as follows:
\begin{equation}
[R^{int}(\theta^{int})]_{r,s} = {\displaystyle \prod_{{j}=1}^{m}{  \exp \left(-   \theta^{int}_{j}  \abs{z^r_{j} - z^s_{j}}^p \right) }}.
\end{equation}
In a similar fashion, constructing $R^{int}(\theta^{int})$ would require the estimation of $m$ non-negative hyperparameters, $\textit{i.e.}$, $\theta^{int} \in \mathbb{R}^m_+$.

\subsection{Correlation matrices for categorical inputs}
\label{subsec:mi_kriging}

For categorical inputs, different choices can be made to build the correlation matrix $R^{cat}(\theta^{cat})$. Some choices are sophisticated and can therefore lead to better GP models, but are known to be computationally expensive (particularly as the number of categorical inputs increases)~\cite{Roustant,Pelamatti}. On the contrary, simple extensions of the well-known continuous kernels based on the Gower distance~\cite{Gower} or on the continuous relaxation techniques~\cite{one-hot} would be less expensive.
In the rest of this section, we will describe three known techniques to build correlation matrices for categorical inputs that are based on kernels. 

\subsubsection{Gower distance based kernel}

The Gower distance based kernel dedicates one hyperparameter per categorical input variable~\cite{Gower,RaulAIAA}. 
Namely, for two given inputs $c^r \in \mathbb{F}^l$ and $c^s \in \mathbb{F}^l$, the Hamming distance, or score, $s$ between the $i^{th}$ component of $c^r$ and $c^s$ is defined as: $s(c_i^r, c_i^s)=0$ if $c_i^r = c_i^s$, otherwise $s(c_i^r, c_i^s)=1$. Thanks to the Hamming distance, one can straightforwardly uses a continuous kernel  to define $R^{cat}(\theta^{cat})$. For instance, in the case of an exponential kernel, the Gower distance based correlation matrix will be given by 
$$
[R^{cat}(\theta^{cat})]_{r,s} =k^{cat}(c^r,c^s,\theta^{cat}) = {\displaystyle \prod_{i=1}^{l}{   \exp \left(- \theta_i^{cat} s(c_i^r,c_i^s)^p  \right) }}.  
$$
Similarly to the continuous and integer correlation matrices, the construction of the categorical correlation matrix based on the Gower distance kernel requires the estimation of $l$ hyperparameters ($\theta^{cat} \in \mathbb{R}^l_+$). Note that, as the Hamming distance can only take the values $0$ and $1$, all the exponential kernels lead to the same result independently of the value of $p$.

\subsubsection{Continuous relaxation based kernel}

To handle categorical variables through continuous relaxation, the design space $ \mathbb{F}^l $ is relaxed to a continuous space $ \Omega^l$  constructed in the following way. 
For a given $i\in \{1,\ldots,l\}$, let $c_i$ be the $i^{th}$ categorical variable with $L_i$ levels, and, for a given input point $c^r$, let $\ell_r^i$ be the index of the level taken by $c^r$ on the variable $i$.
Denote $e_{c^r_i}$ the one-hot encoding~\cite{one-hot} of $c^r_i$ that takes value $0$ everywhere but on the dimension $\ell_r^i$:  $e_{c^r_i} \in \mathbb{R}^{L_i} $ such that $\left( e_{c^r_i} \right)_{\ell^i_r}=1$ and $\left( e_{c^r_i} \right)_k=0$ for $k \neq \ell^i_r$. 
For example, if the $i^{th}$ component is the color with $L_i=3$ levels being $ \{ \mbox{red}, \mbox{blue}, \mbox{green} \}$ and if the $r^{th}$ variable takes value $\mbox{blue}$ ($c_i^r = \mbox{blue}$), then the corresponding index is $\ell_r^i =2$ and the corresponding one-hot encoding is $e_{c^r_i} = (0,1,0)$.
The continuous relaxation idea is as follows. At the beginning we set $ \Omega^l$ to be empty, then, for each $i  \in \{1, \ldots,l\}$, a relaxed one-hot encoding is used for $c_i$. The latter increases the dimension of the relaxed continuous space $ \Omega^l$ by $L_i$ and, at the end of the relaxation, we get the final continuous design space $\Omega^l\subseteq \{0,1\}^{n^l}$, where $n^l=\sum_{i=1}^l L_i >l$. 
Like the Gower distance based kernel, the continuous relaxation based kernel adapts continuous kernels to handle categorical variables, $\textit{i.e.}$, for a couple of categorical inputs $c^r$ and $c^s$,
\begin{equation*}
[R^{cat}(\theta^{cat})]_{r,s} =k^{cat}(c^r,c^s,\theta^{cat}) ={\displaystyle \prod_{i=1}^{l}{ k^{cat}( c_i^r, c_i^s,\theta^{cat})}}= {\displaystyle \prod_{i=1}^{l}{ {\displaystyle \prod_{j=1}^{L_i} k^{cont}( [e_{c_i^r}]_j, [e_{c_i^s}]_{j},\theta^{cat})}}}.    
\end{equation*}
Typically, for an exponential continuous kernel, one has
\begin{equation}
\begin{split}
[R^{cat}(\theta^{cat})]_{r,s} &= {\displaystyle \mathlarger{\prod}_{i=1}^{l} {\displaystyle \mathlarger{\prod}_{j=1}^{L_i} \exp\left(-\theta_{ \sum_{i'=1}^{i-1} {L_{i'}+j}}^{cat}\abs{[e_{c^r_i}]_j - [e_{c^s_i}]_j} ^p\right) }},
\end{split}
\label{eq:kernel}
\end{equation}
and, by using the one-hot encoding structure of $e_{c^r_i}$ and $e_{c^s_i}$, it leads to
\begin{equation*}
\begin{split}
[R^{cat}(\theta^{cat})]_{r,s} &= {\displaystyle \mathlarger{\prod}_{i=1}^{l}  \exp\left(- \theta_{ \sum_{i'=1}^{i-1} {L_{i'}+ \ell^i_r}}^{cat} - \theta_{ \sum_{i'=1}^{i-1} {L_{i'}+ \ell^i_s}}^{cat} \right). } 
\end{split}
\label{eq:kernel:2}
\end{equation*}
Hence, this kernel relies on $n^l=\sum_{i=1}^l L_i$ hyperparameters ($\theta^{cat} \in \mathbb{R}^{n^l}_+$) which can be much more higher than the number of hyperparameters required to build the Gower distance based kernel. Due to one-hot encoding strategy, the value of $p$ is also irrelevant for the construction of the continuous relaxation based kernel.

\subsubsection{Homoscedastic hypersphere kernel}
\label{subsec:ho_hs}
The idea of the homoscedastic hypershere kernel~\cite{Roustant,Pelamatti} is to directly model the correlation matrix instead of looking for a kernel function.
The use of a kernel function guarantees the related correlation matrix $R^{cat}$ to be symmetric positive definite (SPD). However, with the homoscedastic hypershere kernel, one will directly construct an SPD matrix with the desired properties. Namely, for a given $i \in \{1, \ldots, l\}$, let $c^r_{i} $  and $c^s_{i} $  be a couple of categorical variables taking respectively the $\ell^i_r$ and the $\ell^i_s$ level on the categorical variable $c_i$, $[R^{cat}(\theta^{cat})]_{r,s}$ can be formulated in a  level-wise form~\cite{Pelamatti} as: 
\begin{equation}
[R^{cat}(\theta^{cat})]_{r,s}=k^{cat}(c^r,c^s,\theta^{cat}) 
= {\displaystyle \prod_{i=1}^{l}  [R_i(\Theta_i)]_{\ell^i_r,\ell^i_s} }= {\displaystyle \prod_{i=1}^{l}  [C(\Theta_i)C(\Theta_i)^{\top}]_{\ell^i_r,\ell^i_s} }.
\label{eq:homo_HS}
\end{equation}
For all $i \in \{1, \ldots, l\}$, the matrix $C(\Theta_i)\in \mathbb{R}^{L_i \times L_i}$ is lower triangular and built using a hypersphere decomposition~\cite{HS,HS_Jacobi} from a symmetric matrix $\Theta_i \in \mathbb{R}^{L_i \times L_i}$ of hyperparameters. For any $k, k' \in \{1,\ldots, L_i \}$, the  matrix $C(\Theta_i)$ is given by: 
\begin{equation}
 \begin{cases}
 [C(\Theta_i)]_{1,1} = 1, \\
 [C(\Theta_i)]_{k,1} =  \cos\left([\Theta_i]_{k,1}\right) ~\mbox{for any}~  2 \leq k \leq L_i \\
 [C(\Theta_i)]_{k,k'} =  \cos\left([\Theta_i]_{k,k'}\right)  \prod_{j=1}^{k'-1} \sin\left([\Theta_i]_{k,j}\right),~\mbox{for any}~ 2 \le k'< k \leq L_i \\
 [C(\Theta_i)]_{k,k} =    \prod_{j=1}^{k-1} \sin\left([\Theta_i]_{k,j}\right),~\mbox{for any}~ 2 \le k \leq L_i,  \\
\end{cases} 
\label{eq:hy_decomp}
\end{equation}
where the hyperparameters are set such that  $[\Theta_i]_{k,k'} \in [0,\pi]$ for all $1 \le k'< k \leq L_i $. 
For this kernel, the hyperparameters $\theta^{cat}$ can be seen as a concatenation of the set of symmetric matrices, \textit{i.e.}, $\theta^{cat} = \{  \Theta_1, \Theta_2, \ldots, \Theta_l \} $. The construction of this kernel is thus relying on the estimation of $\sum_{i=1}^l \frac{1}{2} L_i (L_i-1) $ hyperparameters. Unlike the previous kernels where the elements of the correlation matrix are non-negative, the correlation values for the homoscedastic hypersphere kernel can be negative, \textit{i.e.}, $[R^{cat}(\theta^{cat})]_{r,s} \in [-1,1]$.

\section{An exponential kernel-based model for categorical inputs}
\label{sec:model_uni}

In this section, we propose an  extension of the classical exponential kernels (used for continuous inputs) to handle categorical variables.
Thanks to the one-hot encoding, we can replace the distance-based approach by an hyperparameter-based approach. 
This extension will naturally lead to a generalization of both continuous relaxation and Gower distance based kernels.  

Distance based approaches (like Gower distance or continuous relaxation) can not model every possible correlation between the various categorical choices. Therefore, these methods do not lead to an exhaustive GP model but to an imprecise approximation. In what follows, we propose to introduce a new formulation that includes a correlation matrix so that we could reach a higher accuracy for the resulting distance-based GP model.
To begin with, the continuous relaxation kernel described in~\eqnref{eq:kernel} can be reformulated as:
\begin{equation}
\begin{split}
[R^{cat}(\theta^{cat})]_{r,s}
&= {\displaystyle \mathlarger{\prod}_{i=1}^{l} {\displaystyle \mathlarger{\prod}_{j=1}^{L_i} \exp\left(- \abs{[e_{c^r_i} - e_{c^s_i}]_j} ^{p/2} [{\Theta}_{i}]_{j,j} \abs{[e_{c^r_i} - e_{c^s_i}]_j} ^{p/2}\right) }},
\end{split}
\label{eq:CR2}
\end{equation}
where, for all $i=1, \ldots,l$, the matrix $\Theta_i \in \mathbb{R}^{L_i\times L_i}$ is diagonal  such that $[\Theta_i]_{j,j}=\theta_{ \sum_{i'=1}^{i-1} {L_{i'}+j}}^{cat} \in \mathbb{R}_+$, and $\theta^{cat}$ is defined as the list of hyperparameter matrices $\theta^{cat} = \{\Theta_1, \ldots, \Theta_l\}$.
The idea of the new kernel is the following: we start from the reformulation of the continuous relaxation kernel
of~\eqnref{eq:CR2}. Then, as for the kernel of~\eqnref{eq:homo_HS}, we consider, for every categorical variable $i=1, \ldots,l$, a SPD matrix $\Phi( \Theta_i) \in  \mathbb{R}^{L_i\times L_i}$ used to build a kernel associated with the correlation matrix $R_i(\Phi(\Theta_i))$.
Let $c^r_{i} $  and $c^s_{i} $  be a couple of categorical variables taking respectively the $\ell^i_r$ and the $\ell^i_s$ level of the variable $c_i$, we set
\begin{equation}
[R^{cat}(\theta^{cat})]_{r,s} = {\displaystyle \prod_{i=1}^{l}  [R_i(\Phi(\Theta_i))]_{\ell^i_r,\ell^i_s} },
\label{eq:nat_ext_gaussian_ker:0}
\end{equation}
and, for all $i=1, \ldots,l$, one has
\begin{equation}
[R_i(\Phi(\Theta_i))]_{\ell^i_r,\ell^i_s}  
=   {\displaystyle \prod_{j=1}^{L_i}  {\displaystyle \prod_{j'=1}^{L_i}  
\exp \left(-   \abs{ [e_{c^r_i} -  e_{c^s_i}]_j}^{p/2} [\Phi(\Theta_i)]_{j,j'}  \abs{ [e_{c^r_i}   -  e_{c^s_i}]_{j'}}^{p/2} \right) }},
\label{eq:nat_ext_gaussian_ker}
\end{equation}
where $\ell_r^i$ and $\ell_s^i$ are the indices of the levels taken by the variables $c^r$ and $c^s$, respectively, on the $i^{th}$ categorical variable and the coefficient $[\Phi(\Theta_i)]_{{\ell^i_r},{\ell^i_s}}$ is characterizing the correlation between these two levels. 

\begin{remark}
One can easily see that \eqnref{eq:nat_ext_gaussian_ker} generalizes the continuous relaxation approach. In fact, by setting  $\Phi( \Theta_i)=\Theta_i$ to be a diagonal matrix, we recover~\eqnref{eq:CR2}.
\end{remark}

Now, by using the one-hot encoding nature of the vectors $e_{c_i^r}$ and $e_{c_i^s}$, we get naturally what follows. Namely, if $c^r_i =c^s_i$, one deduces that $ [R_i(\Phi(\Theta_i))]_{\ell^i_r,\ell^i_s}   = \exp (0) = 1$. Otherwise, if $ c^r_i \neq c^{s}_i$, we get
\begin{equation}
\begin{split}
[R_i(\Phi(\Theta_i))]_{\ell^i_r,\ell^i_s}   &= \exp  \left( - {\displaystyle \sum_{j=1}^{L_i} {\displaystyle \sum_{j'=1}^{L_i}{    \abs{  [e_{c^r_i} -  e_{c^s_i}]_j }^{p/2} [\Phi(\Theta_i)]_{j,j'}  \abs{[e_{c^r_i}   -  e_{c^s_i}]_{j'}}^{p/2}  }}} \right) \\
&= \exp \left( - \left([\Phi(\Theta_i)]_{{\ell^i_r},{\ell^i_r}}+[\Phi(\Theta_i)]_{{\ell^i_s},{\ell^i_s}} +[\Phi(\Theta_i)]_{{\ell^i_r},{\ell^i_s}} +[\Phi(\Theta_i)]_{{\ell^i_s},{\ell^i_r}} \right) \right)\\ 
&= \exp \left( - [\Phi(\Theta_i)]_{{\ell^i_r},{\ell^i_r}}-[\Phi(\Theta_i)]_{{\ell^i_s},{\ell^i_s}} - 2[\Phi(\Theta_i)]_{{\ell^i_r},{\ell^i_s}} \right). 
\end{split}
\label{eq:mat_ker_i}
\end{equation}

\begin{remark}
Note that the resulting correlation matrix $R_i(\Phi(\Theta_i))$ does not depend on the chosen parameter $p$ (used within the definition of the exponential kernels). Therefore, in our case, when dealing with categorical variables kernels, there will be no distinction between squared or absolute exponential  kernels. 
\end{remark}

In addition, as far as the matrices $\Theta_i$ respect a specific parameterization, we will show that our approach guarantees that the correlation matrix $R$ is SPD with a unit diagonal and off-diagonal terms values in $[0,1]$~\cite{PDUDE}. In general, the latter properties are required to be satisfied by the correlation matrices. Otherwise, one may get numerical issues to build the GP model, see~\eqnref{eq:mean:GP} and~\eqnref{eq:std:GP}. For that purpose,
for a given $i \in \{1,\ldots, l\}$, we propose to use the following parameterization for the hyperparameter matrix $\Phi(\Theta_i)$:
\begin{equation}
\begin{split}
[\Phi(\Theta_i)]_{j,j} &:= [\Theta_i]_{j,j} ~  \geq 0  \\
[\Phi(\Theta_i)]_{j,j'} &:= \frac{\log \epsilon }{2} ([C(\Theta_i) C(\Theta_i)^\top]_{j,j'} -1)  ~~~\mbox{if}~{j\neq j'},   \\
\end{split}
\label{eq:hs_exp}
\end{equation}
where the parameter $\epsilon$ is chosen as a small positive tolerance ($ 0<\epsilon \ll 1$) and the matrix
$C(\Theta_i)$ is a Cholesky lower triangular matrix that relies on the symmetric matrix $\Theta_i$ of $L_i(L_i-1)/2$ elements in $[0,\frac{\pi}{2}]$. The elements of  $C(\Theta_i)$ represent the coordinates of a point on the surface of a unit radius sphere as in~\cite{ Roustant,Pelamatti}. They are described in~\eqnref{eq:hy_decomp}.  Note that, by taking into consideration the symmetry of the matrix $\Theta_i$, the total number of hyperparameters for the categorical variable $i$ is $\frac{L_i(L_i+1)}{2}$.

In the next theorem, we will show that the parameterization given by \eqnref{eq:hs_exp} guarantees the desirable properties for the correlation matrices $R_i(\Theta_i)$ and therefore for the matrix $R^{cat}$. In particular, we will show  that the matrix $R^{cat}(\theta^{cat})$ is SPD with elements in $[0,1]$, $\textit{i.e.}$, for all $s,r \in \{1,\ldots, n_t\}$,~ $[R^{cat}(\theta^{cat})]_{r,s} \in [0,1]$.

\begin{theorem}
\label{th:definiteness}
 Assume that, for all $i \in \{1, \ldots, l\}$, $ \Phi(\Theta_i)$ satisfies the parameterization of~\eqnref{eq:hs_exp}.
Then the matrix $R^{cat}(\theta^{cat})$, given by \eqnref{eq:nat_ext_gaussian_ker:0}, is SPD with  elements in $[0,1]$.
\end{theorem}
\begin{proof}Indeed, for all $i \in \{1, \ldots, l\}$, by using \eqnref{eq:mat_ker_i} and \eqnref{eq:hs_exp}, one has
\begin{eqnarray*}
[R_i(\Phi(\Theta_i))]_{\ell^i_r,\ell^i_s}&=&[W_i]_{\ell^i_r,\ell^i_s} [T_i]_{\ell^i_r,\ell^i_s}, ~~~~\mbox{if}~ \ell^i_r \neq \ell^i_s \\
[R_i(\Phi(\Theta_i))]_{\ell^i_r,\ell^i_r}&=& 1,
\end{eqnarray*}
where $[W_i]_{\ell^i_r,\ell^i_s}=\exp \left( - [\Theta_i]_{{\ell^i_r},{\ell^i_r}}-[\Theta_i]_{{\ell^i_s},{\ell^i_s}}\right)$ and $[T_i]_{\ell^i_r,\ell^i_s}=\exp \left( - 2 [\Phi(\Theta_i)]_{{\ell^i_r},{\ell^i_s}} \right)$.
The matrix $R_i(\Phi(\Theta_i))$ is thus defined as a Hadamard product  ($\textit{i.e.}$, element-wise product of matrices)~\cite{hadamard}. Hence, by application of the Schur product theorem~\cite[Lemma~3.7.1]{schurmatapp}, it suffices to show that the matrices $W_i$ and $T_i$ are SPD to prove that $R_i$ is also SPD. 
Taking into account that, for all $s,r \in \{1,\ldots, n_t\}$, $e_{c^r_i}$ and $e_{c^s_i}$ are one-hot encoding elements of $\mathbb{R}^{L_i}$,  the matrix $W_i$ corresponds to the correlation matrix associated with the exponential kernel in the continuous space, \textit{i.e.},
$$[W_i]_{\ell^i_r,\ell^i_s}=\exp \left( - [\Theta_i]_{{\ell^i_r},{\ell^i_r}}-[\Theta_i]_{{\ell^i_s},{\ell^i_s}}\right)={\displaystyle \prod_{j=1}^{L_i} {  \exp \left(-   [\Theta_i]_{j,j}  \abs{ [e_{c^r_i}   -  e_{c^s_i}]_{j} }^p \right) }}.$$ 
Hence, since  the diagonal elements of $\Theta_i$ are positive, the matrix $W_i$ is SPD. In fact, the kernel function $\phi(x) = \exp [- \theta |x|^p]$ is positive definite for a given positive $\theta$ if $0<p \leq 2$~\cite[Corollary~3]{Schoenberg}. 
Regarding the matrix $T_i$, by using the fact that $\Theta_i$ satisfies~\eqnref{eq:hs_exp}, one has
\begin{equation}
\label{eq:EHH}
[T_i]_{\ell^i_r,\ell^i_s}= \epsilon  \exp \left( -(\log \epsilon) [C(\Theta_i) C(\Theta_i)^{\top}]_{\ell^i_r,\ell^i_s} \right).
\end{equation}
For an $\epsilon \in (0,1)$,  the matrix $-(\log \epsilon) [C(\Theta_i) C(\Theta_i)^{T}]$ is SPD as a Cholesky like-decomposition matrix. Thus, $T_i$ is also SPD as the Hadamard exponential of an SPD matrix~\cite[Theorem 7.5.9]{horn2012matrix}.

For the second part of the proof, the matrix $C(\Theta_i)$ is constructed by hypersphere decomposition such that the values of $C(\Theta_i) C(\Theta_i)^\top$ belong to $[0,1]$~\cite{hypersphere}. Hence,
\begin{equation*}
\begin{split}
 &[T_i]_{\ell^i_r,\ell^i_s}=\epsilon  \exp \left(  -(\log \epsilon)  [C(\Theta_i) C(\Theta_i)^\top]_{\ell^i_r,\ell^i_s} \right) \geq \epsilon  \exp 0 =  \epsilon,  \\ 
&[T_i]_{\ell^i_r,\ell^i_s}=\epsilon  \exp \left(  -(\log \epsilon)  [C(\Theta_i) C(\Theta_i)^\top]_{\ell^i_r,\ell^i_s} \right) \leq \frac{\epsilon}{\epsilon} =1.
\end{split}
\end{equation*}
Also, the elements of the matrix $W_i$ are in $[0,1]$ since the diagonal elements of $\Theta_i$ are chosen to be positive. Consequently, the extra-diagonal elements of $R_i$ are in $[0,1]$.
Finally, the Hadamard product being conservative for those two latter properties, one concludes that the correlation matrix $R^{cat}(\theta^{cat})$ is SPD and all its elements are in $[0,1]$. 
\end{proof}

\begin{remark}
\label{rmk:2}
For a given small $\epsilon>0$, the transformation
\begin{equation}
\label{eq:bijection_alpha}
\alpha \to \epsilon \exp[-\log(\epsilon) (\alpha-1)]    
\end{equation}
is a bijection over $[\epsilon,1]$, thus one can deduce that (when $[T_i]_{j,j'}> \epsilon$ for all $j,j'$) there exists a unique  matrix $\hat \Theta_i$ such that 
$$
T_i = C(\hat \Theta_i)C(\hat \Theta_i)^{\top}.
$$
This, in particular, shows that if we set $W_i$ to identity in our parameterization, then, as far as the correlations are larger than $\epsilon$, the homoscedastic hypersphere parameterization of Zhou \textit{et al.}~\cite{Zhou,Pelamatti}  is equivalent to our proposed one.
\end{remark}

In the next theorem, using the hypersphere decomposition properties~\cite{HS_Jacobi}, we will show that the correlation matrix $R_i$, as given by~\eqnref{eq:mat_ker_i}, can be built in an equivalent way without the diagonal elements of the matrix  $\Phi(\Theta_i)$. Such result is of high interest as it reduces the number of hyperparameters from $\frac{L_i(L_i+1)}{2}$ to  $\frac{L_i(L_i-1)}{2}$ per categorical variable $i$ without any loss in the accuracy in the final model.
\begin{theorem}
\label{lemma_eq_ho_hs}
The correlation matrix $R_i$, as given by~\eqnref{eq:mat_ker_i}, can be rewritten as follows
\begin{equation}
\label{new:R}
\begin{split}
[R_i(\Phi(\bar \Theta_i))]_{\ell^i_r,\ell^i_s}=& \exp \left( - 2 [\Phi(\bar \Theta_i)]_{{\ell^i_r},{\ell^i_s}} \right), ~~~~\mbox{if}~ \ell^i_r \neq \ell^i_s \\
[R_i(\Phi(\bar \Theta_i))]_{\ell^i_r,\ell^i_r}=& 1,
\end{split}
\end{equation}
where $[\Phi(\bar \Theta_i)]_{{\ell^i_r},{\ell^i_s}}= \frac{\log \epsilon }{2} ([C(\bar\Theta_i) C(\bar\Theta_i)^\top]_{{\ell^i_r},{\ell^i_s}} -1)$ and $\bar\Theta_i$ is a symmetric matrix whose diagonal elements are set to zero (\textit{i.e.}, $[\bar\Theta_i]_{j,j} =0$ for all $j=1, \ldots, L_i$).
\end{theorem}
\begin{proof}
Indeed, by using the hypersphere decomposition~\cite{HS_Jacobi}, any SPD matrix $T_i(\Theta_i)$ with unitary diagonal and values in $[\epsilon,1]$ can be modeled as $T_i(\Theta_i) = [C({\hat \Theta}_i) C({\hat\Theta}_i)^\top]$ from a certain symmetric matrix ${\hat \Theta_i}$ without using additional diagonal elements (\textit{i.e.}, $[\hat \Theta_i]_{j,j}=0$ for all $j=1, \ldots, L_i$).
Thus, using the fact that $R_i(\Theta_i)$ is written as the image of this SPD matrix $T_i(\Theta_i)$ by the element-wise transformation of ~\eqnref{eq:bijection_alpha} bijective over $[\epsilon,1]$, one can deduce that there must exist a symmetric matrix ${\bar \Theta_i}$ whose diagonal elements are set to zero (\textit{i.e.}, $[\bar \Theta_i]_{j,j}=0$ for all $j=1, \ldots, L_i$) such that 
\begin{equation*}
\begin{split}
[R_i(\Phi(\bar \Theta_i))]_{\ell^i_r,\ell^i_s}=& \exp \left( - 2 [\Phi(\bar \Theta_i)]_{{\ell^i_r},{\ell^i_s}} \right), ~~~~\mbox{if}~ \ell^i_r \neq \ell^i_s \\
[R_i(\Phi(\bar \Theta_i))]_{\ell^i_r,\ell^i_r}=& 1,
\end{split}
\end{equation*}
where $[\Phi(\bar \Theta_i)]_{{\ell^i_r},{\ell^i_s}}= \frac{\log \epsilon }{2} ([C(\bar\Theta_i) C(\bar\Theta_i)^\top]_{{\ell^i_r},{\ell^i_s}} -1)$. 
\end{proof}
 
 In what comes next, we will refer to our kernel when it uses the parameterization of~\eqnref{new:R} as the Exponential Homoscedastic Hypersphere (EHH) kernel (see Remark~\ref{rmk:2}). We will call the original parameterization of the correlation matrix (as given by~\eqnref{eq:mat_ker_i}) as the Fully Exponential (FE) kernel. Note that, as explained in~\ref{apendix:EHH2CR}, whenever the matrix $\Theta_i$ is chosen to be diagonal, the matrix $\Phi(\Theta_i)$ used within FE kernel will also be diagonal. Thus, we are able to recover the continuous relaxation kernel~\cite{GMHL} and  this parameterization will be called the Continuous Relaxation (CR) kernel. Similarly, if we choose $\Theta_i $ to be of the form $\theta_i \times I_{L_i}$ where $\theta_i \in \mathbb{R}^+$ and $I_{L_i}$ is the identity matrix of size $L_i$, we are able to recover the Gower distance based kernel~\cite{Gower}: this parameterization will be called the Gower Distance (GD) kernel. 
 
 As mentioned earlier, the EHH kernel is similar to the FE kernel. Therefore, we can deduce that the EHH kernel generalizes the CR kernel and also that the CR kernel generalizes the GD kernel. Table~\ref{tab:General scheme} gives all the details associated with the four categorical kernels described above, \textit{i.e.}, GD, CR, EHH and FE.

\begin{table}[H]
   \caption{Description of the four categorical kernels (GD, CR, EHH and FE) using our proposed exponential parameterization.}
   \vspace*{-0.2cm}
   \begin{center}
   \resizebox{1\columnwidth}{!}{%
      \begin{tabular}{|c|c|c|c|}
       \hline
  Kernel &    $\Theta_i=$ & 

 $[R_i(\Theta_i)]_{{\ell^i_r},{\ell^i_s}}=$ 
          
       &  \# of Hyperparam. \\
       \hline 
 & & &  \\
 GD 
 & 
 \scalebox{1.2}{$ \displaystyle \frac{\theta_i}{2} \  \times $  } $  
\begin{bmatrix}
1 & \textcolor{white}{9} & \hspace{2em} { \textbf{\textit{ Sym.}}}  \textcolor{white}{9} & \\
0  &1 & \textcolor{white}{9} \\
\vdots & \ddots &  \ddots & \textcolor{white}{9}  \\
0 &  \ldots & 0 &1 \\
\end{bmatrix}$  &
$ \exp \left( - [\Phi(\Theta_i)]_{{\ell^i_r},{\ell^i_r}}-[\Phi(\Theta_i)]_{{\ell^i_s},{\ell^i_s}}  \right) $
&1\\
 & & &  \\
\hline 
 & & &  \\
 CR 
 & 
  $   
\begin{bmatrix}
[\Theta_i]_{1,1}  & \textcolor{white}{9} & \textcolor{white}{9} & \hspace{-2em} { \textbf{\textit{ Sym.}}} \\
0 &[\Theta_i]_{2,2}  & \textcolor{white}{9} \\
\vdots & \ddots &  \ddots & \textcolor{white}{9}  \\
0 &  \ldots & 0 & [\Theta_i]_{L_i,L_i} \\
\end{bmatrix}$ & 
  $ \exp \left( - [\Phi(\Theta_i)]_{{\ell^i_r},{\ell^i_r}}-[\Phi(\Theta_i)]_{{\ell^i_s},{\ell^i_s}}  \right) $
& $ L_i$ \\
 & & &  \\
\hline 
 & & &  \\
 EHH  & 
$ 
\begin{bmatrix}
0 & \textcolor{white}{9} & \hspace{2em} { \textbf{\textit{ Sym.}}}  \textcolor{white}{9} & \\
[\Theta_i]_{1,2}  &0 & \textcolor{white}{9} \\
\vdots& \ddots &  \ddots & \textcolor{white}{9}  \\
[\Theta_i]_{1,L_i} &  \ldots & [\Theta_i]_{L_i-1,L_i} &0 \\
\end{bmatrix}$ &  
$\exp \left(  - 2  [\Phi(\Theta_i)]_{{\ell^i_r},{\ell^i_s}} \right)$
&  $\frac{1}{2} L_i(L_i-1) $ \\
 & & &  \\
     \hline
      & & &  \\
FE & 

$\begin{bmatrix}
[\Theta_i]_{1,1} & \textcolor{white}{9} & \hspace{2em} { \textbf{\textit{ Sym.}}}  \textcolor{white}{9} & \\
[\Theta_i]_{1,2}  & [\Theta_i]_{2,2} & \textcolor{white}{9} \\
\vdots &\ddots & \ddots & \textcolor{white}{9}  \\
[\Theta_i]_{1,L_i} &  \ldots & [\Theta_i]_{L_i-1,L_i} &[\Theta_i]_{L_i,L_i} \\ 
\end{bmatrix} $  & 
$\exp \ ( - [\Phi(\Theta_i)]_{{\ell^i_r},{\ell^i_r}}-[\Phi(\Theta_i)]_{{\ell^i_s},{\ell^i_s}}- 2[\Phi(\Theta_i)]_{{\ell^i_r},{\ell^i_s}} )$ &

  $\frac{1}{2} L_i(L_i+1) $ \\
   & & &  \\
        \hline  
      \end{tabular}
      }
   \end{center}
   \label{tab:General scheme}
\end{table}

{\color{black}
To sum up, we have seen that the HH kernel can be more general than the EHH  one (as it can deal with negative correlations) and that EHH generalizes both CR and GD. All these categorical models can be unified in a single formulation as follows. 
For each $i \in \{1, \ldots, l\}$,  a hyperparameter matrix $\Theta_i$ is associated with each variable $c_i$, i.e.,
$$\Theta_i= \begin{bmatrix}
[\Theta_i]_{1,1} & \textcolor{white}{9} & \hspace{2em} { \textbf{\textit{ Sym.}}}  \textcolor{white}{9} & \\
[\Theta_i]_{1,2}  & [\Theta_i]_{2,2} & \textcolor{white}{9} \\
\vdots &\ddots & \ddots & \textcolor{white}{9}  \\
[\Theta_i]_{1,L_i} &  \ldots & [\Theta_i]_{L_i-1,L_i} &[\Theta_i]_{L_i,L_i} \\ 
\end{bmatrix}.$$
The correlation term $[R_i(\Theta_i)]_{\ell^i_r,\ell^i_s} $ associated with $c_i$ can be formulated in the following level-wise form: 
\begin{eqnarray*}
[R_i(\Theta_i)]_{\ell^i_r,\ell^i_s} &= &\   \kappa ( 2[ \Phi(\Theta_i) ]_{{ \ell_i^r},{\ell_i^s}} ) \  \kappa ( [ \Phi(\Theta_i) ]_{{ \ell_i^r},{\ell_i^r}} ) \  \kappa ( [ \Phi(\Theta_i) ]_{{ \ell_i^s},{\ell_i^s}} ), 
\label{eq:homogeneous}
\end{eqnarray*}

where $\kappa$ can be set either to any positive definite kernel (to get  GD, CR, EHH or FE kernels) or to identity (to get HH kernel). The transformation 
function $\Phi(.)$ is selected such that, for any SPD matrix $\Theta_i$, the output matrix $\Phi(\Theta_i)$ is also SPD. \tabref{tab:kernels} gives a list of possible choices for $\Phi$ when the is function $\kappa$ is set to exponential or identity. 
For all categorical variables $i \in \{1, \ldots, l\}$, the matrix $C(\Theta_i)\in \mathbb{R}^{L_i \times L_i}$ (lower triangular) is  built using a hypersphere decomposition. 

\begin{table}[H]
\caption{Kernels using different choices for the function $\Phi$.}
\begin{tabular}{c|c|l}
Kernel & $\kappa(\phi)$  &  \hspace{3cm} ${\centering \Phi(\Theta_i)}$  \\
\hline
{GD}   &  $\exp(-\phi) $ & ${ \displaystyle [\Phi(\Theta_i)]_{j,j} := \frac{1}{2}  \theta_{i} \quad  ~;~ [\Phi(\Theta_i)]_{j \neq j'} := 0 }$    \\
{CR}  & $\exp(-\phi) $ &  $  { \displaystyle [\Phi(\Theta_i)]_{j,j} := [\Theta_i]_{j,j} ~;~ [\Phi(\Theta_i)]_{j \neq j'} := 0 } $   \\
{EHH}  & $\exp(-\phi)$ & 
$  { \displaystyle [\Phi(\Theta_i)]_{j,j} := 0 \quad \quad ~;~ [\Phi(\Theta_i)]_{j \neq j'} := \frac{\log \epsilon }{2} ([C(\Theta_i) C(\Theta_i) ^\top]_{j,j'} -1)  }$ \\
{HH} &  $\phi$ &    $  { \displaystyle [\Phi(\Theta_i)]_{j,j} := 1 \quad \quad ~;~ [\Phi(\Theta_i)]_{j \neq j'} :=  \frac{1}{2} [C(\Theta_i) C(\Theta_i)^\top]_{j,j'} }$  \\
\end{tabular}
\label{tab:kernels}
\end{table}
}

In the next section, we will see how these kernels perform on different test cases. In particular, we study numerically the trade-off between the kernel efficiencies and their respective computational efforts (related directly to the number of hyperparamters).

\section{Results and discussion}
\label{sec:Results}

In this section, we propose several illustrations and comparisons on \textcolor{black}{three} different test cases (from  \textcolor{black}{2} to 10 continuous variables and 1 or 2 categorical variables up to 12 levels) to show the interest of our method and the equivalence with other kernels from the literature.
The likelihood value and the approximate errors are the quantities of interest considered to compare different correlation kernels.
\subsection{Implementation details}

The optimization of the likelihood as a function of the hyperparameters needs a performing gradient-free algorithm, in this work, we are using COBYLA~\cite{COBYLA} to maximize this quantity from the Python library Scipy with default termination criterion related to the trust region size. As COBYLA is a local search algorithm, a multi-start technique is used. Our models and their implementation are available in the toolbox SMT v2.0\footnote{\url{https://smt.readthedocs.io/en/latest/}}~\cite{SMT2019}.  By default, in SMT, the number of starting points for COBYLA is equal to 10 with evenly spaced starting points.

A simple noiseless Kriging with a constant prior model for the GP is used. We recall that the absolute exponential kernel and the squared exponential kernel are similar for categorical variables and  differ only for the continuous ones. The correlation values range between 2.06e-9 and 0.999999 for both continuous and categorical hyperparameters. Therefore, the constant $\epsilon$ is chosen to correspond to a correlation value of 2.06e-9.
The random DoEs are drawn by Latin Hypercube Sampling (LHS)~\cite{LHS} and  the validation sets are given by some evenly spaced points. 

\subsection{ \textcolor{black}{Analytic validation: a categorical cosine problem ($n= 1$, $m=0$, $l=1$ and $L_1=13$)}}

In this section,  we consider the categorical cosine problem, from~\cite{Roustant}, to illustrate the behaviour of our proposed kernels. In this problem, the objective function $f$ depends on a continuous variable in $[0,1]$ and on a categorical variable with 13 levels.~\ref{subsec:cosine} provides a detailed description of this function. Let $w= (x,c)$ be a given point with  $x$ being the continuous variable and $c$ being the categorical variable, $c \in \{1, \ldots, 13\}$. There are two groups of curves corresponding to levels 1 to 9 and levels 10 to 13 with strong positive within-group correlations, and strong negative between-group correlations. 

\begin{figure}[h]
\begin{center}

	\subfloat[GD kernel]{
      \centering 
		\includegraphics[clip=true, height=4.5cm, width=7.2cm]{  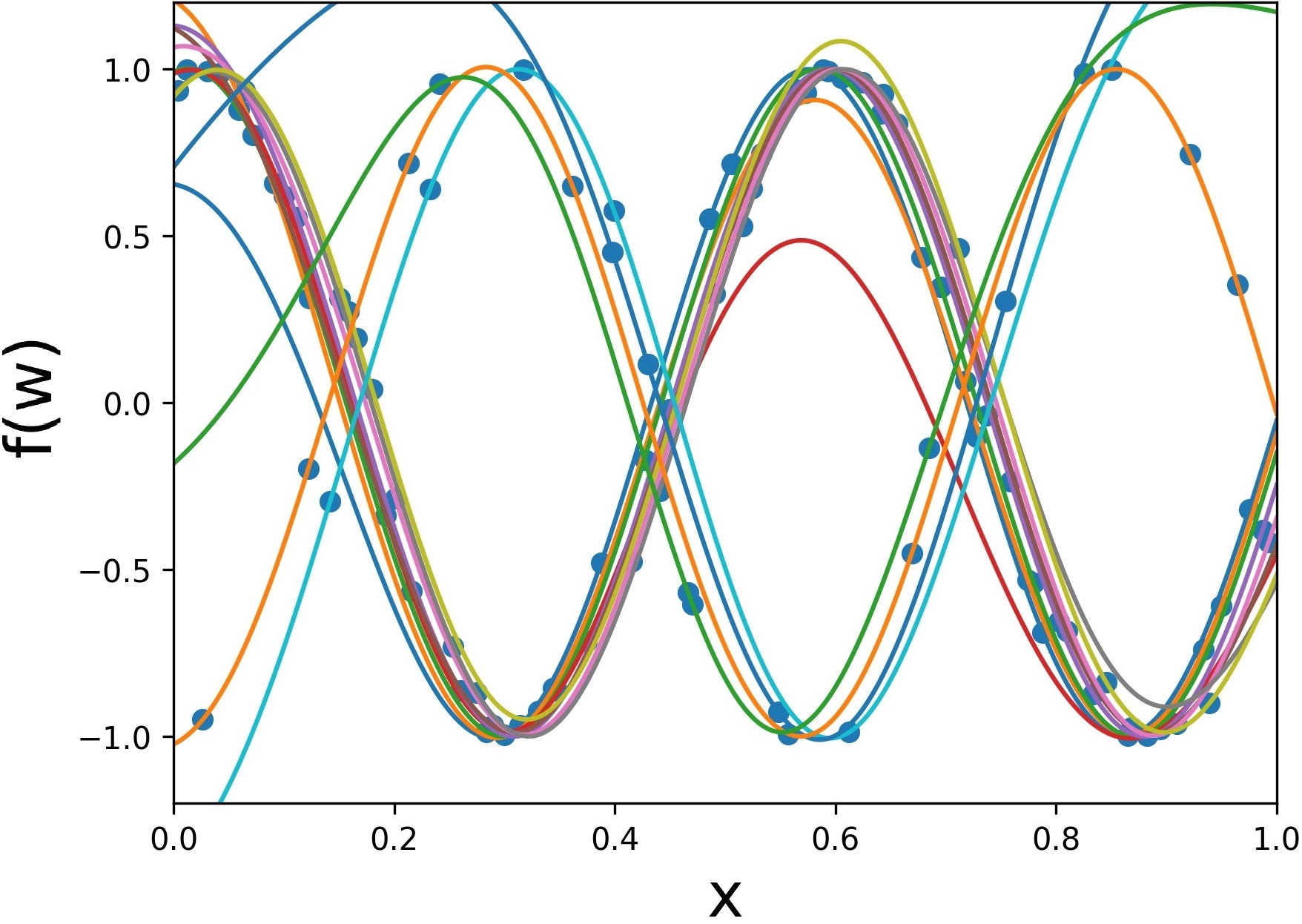}
     }   
        \subfloat[CR kernel]{
      \centering
	\includegraphics[clip=true, height=4.5cm, width=8cm]{  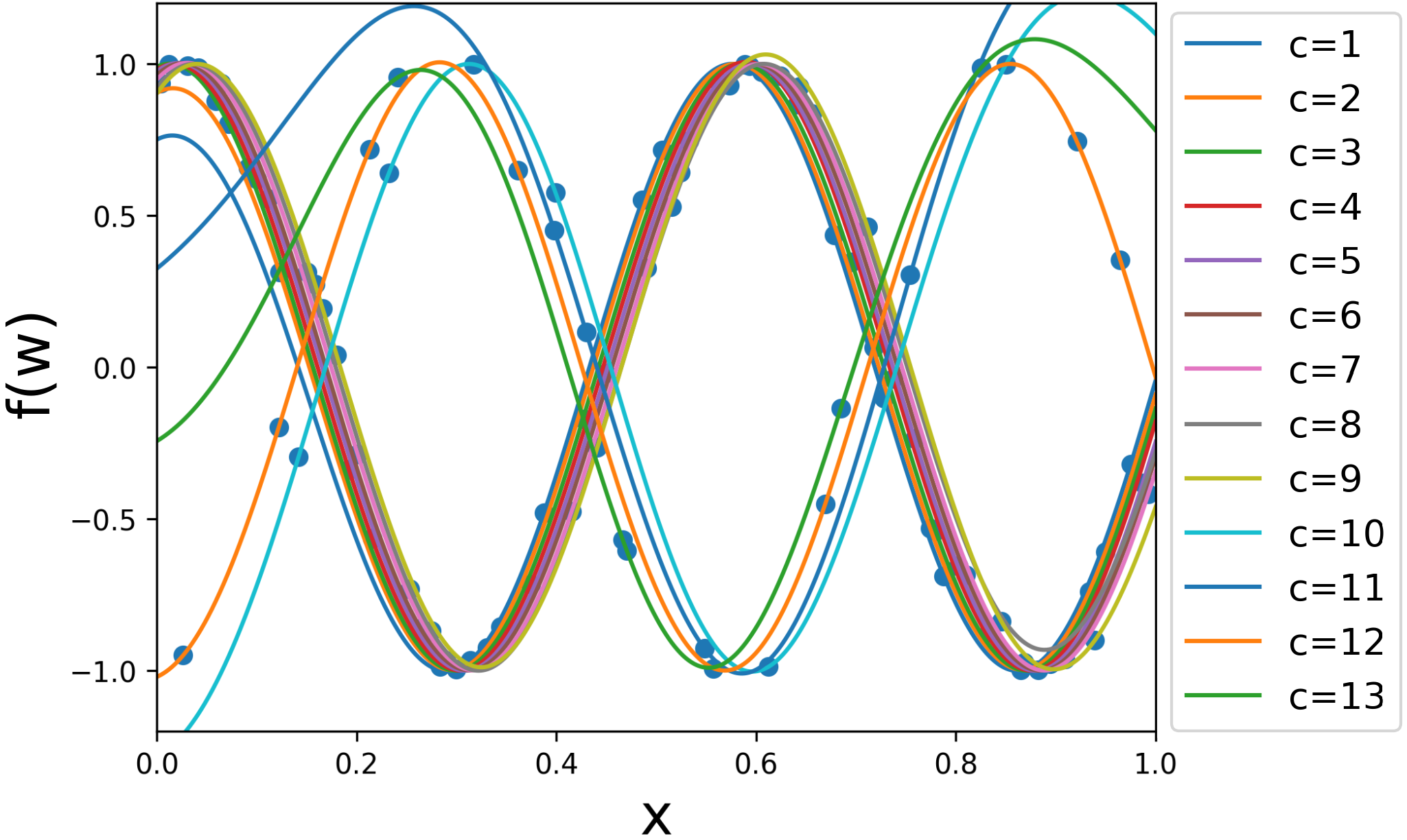}} 
	
        \subfloat[EHH kernel]{
      \centering
	\includegraphics[clip=true, height=4.5cm, width=7.2cm]{  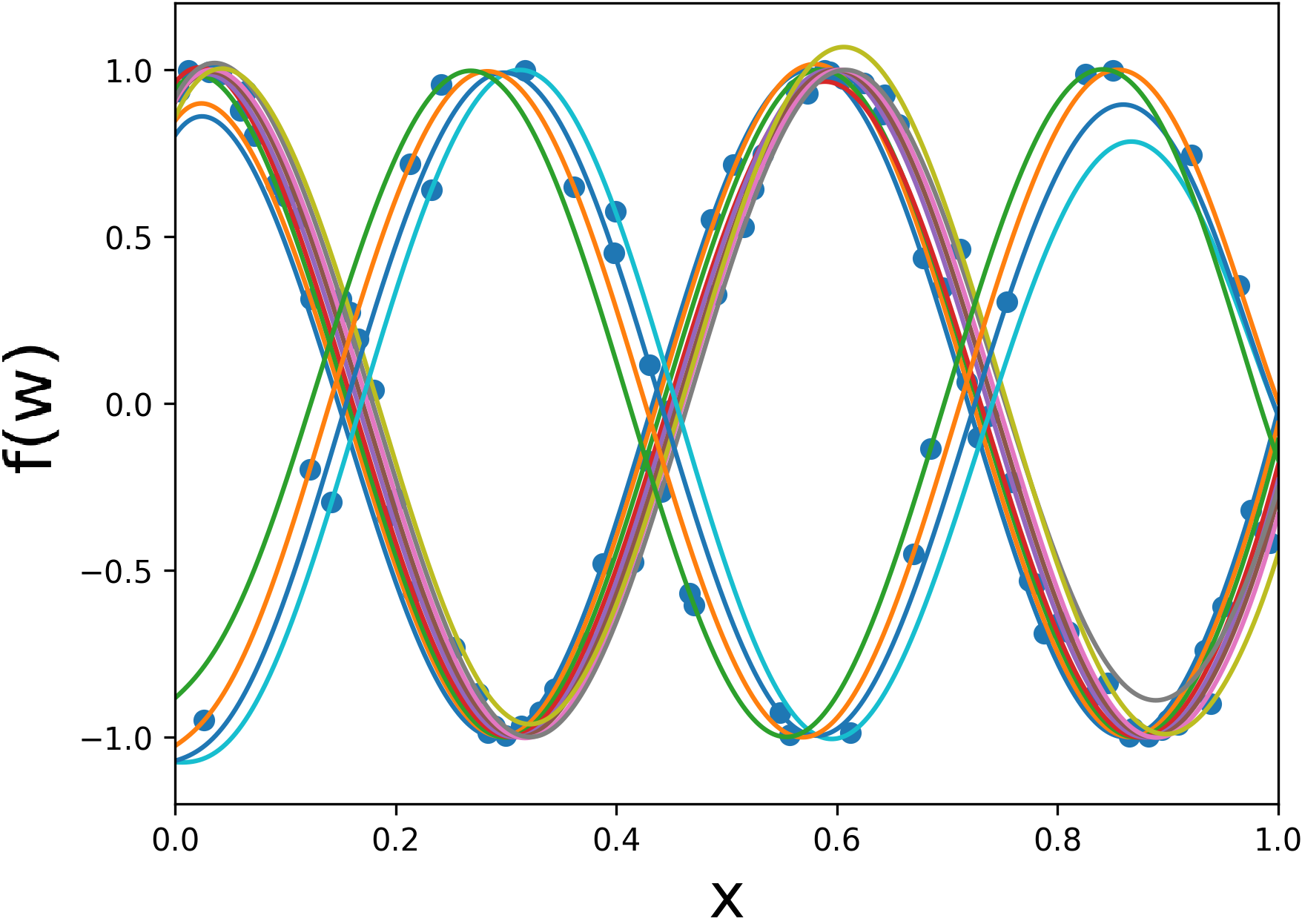}} 
 \subfloat[HH kernel]{
  \centering
	\includegraphics[clip=true, height=4.5cm, width=7.2cm]{  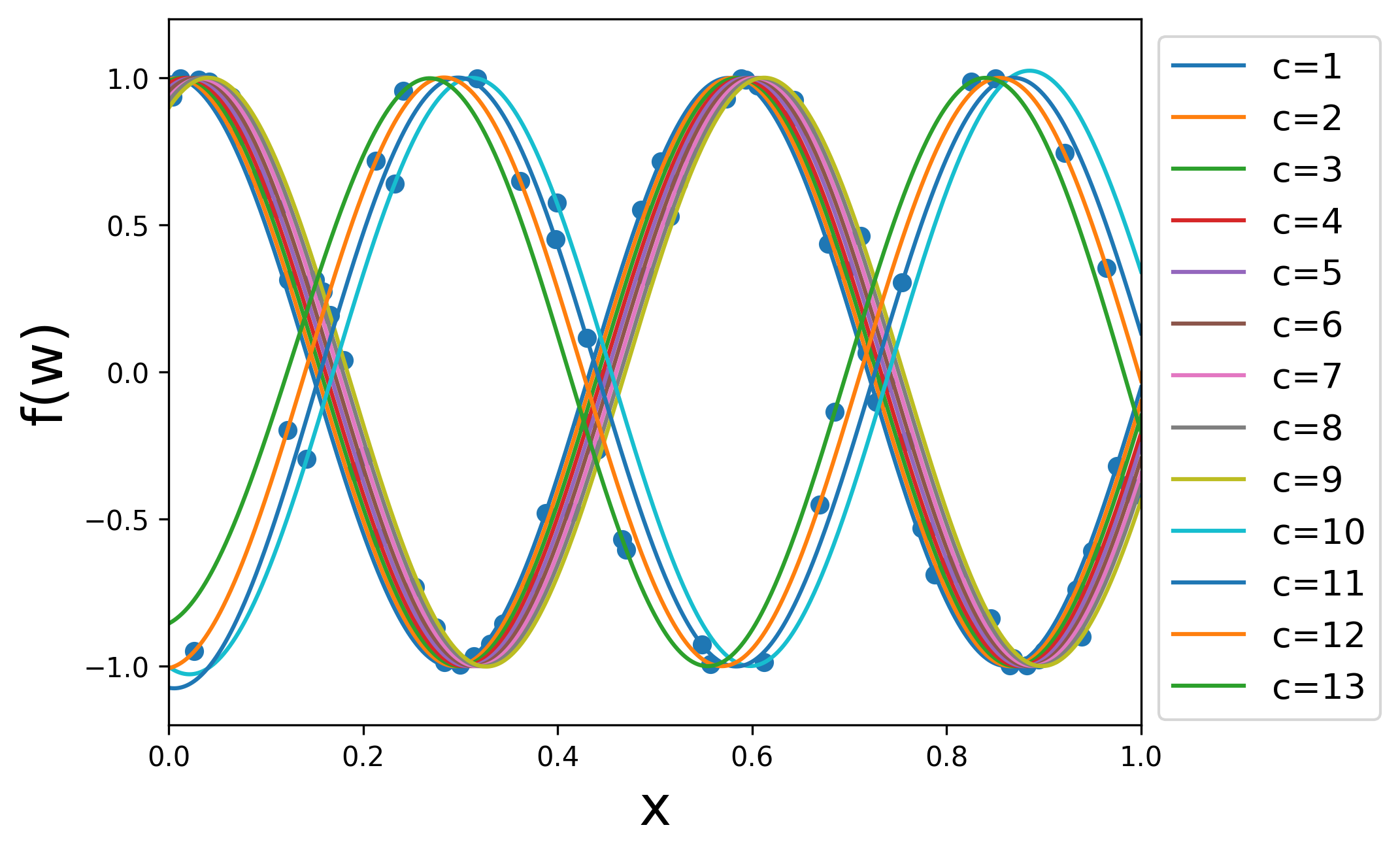} 
    \label{Roustant_hs_gsh_theo}
     }
\caption{Mean predictions for our proposed model using different kernels over the matrix $\Theta_1$ for the cosine problem with a 98 point DoE. }
\label{models_Roustant}
\end{center} 
\end{figure}

\textcolor{black}{
In this example, the number of relaxed dimensions for continuous relaxation is 14. A LHS DoE with 98 points ($14\times 7$, if 7 points per dimension are considered) is chosen to built the Gaussian process models. 
The associated mean posterior models are shown on~\figref{models_Roustant} for GD, CR, EHH and HH. 
 The number of hyperparameters to optimize is therefore $2$ for GD, $14$ for CR and $79$ for EHH and HH as indicated in~\tabref{tab:resRoustant}. 
 }
{\color{black}\figref{models_Roustant} shows that the predicted values remain properly within the interval $[-1,1]$ only with HH and EHH kernels. Therefore, these kernels seem to be better modeling methods. 

To better assess the accuracy of each kernel, we compute the root mean square error (RMSE) and the predictive variance adequacy (PVA)~\cite{PVA} are respectively given by 
\begin{equation*}
\mbox{RMSE}= \sqrt{\underset{i=1}{\overset{n}{\sum}} \frac{1}{n} \left( \hat{f}(w_i) - f(w_i) \right)^2} \quad \mbox{and} \quad
\mbox{PVA}= \log \left( {\underset{i=1}{\overset{n}{\sum}} \frac{1}{n} \frac{\left( \hat{f}(w_i) - f(w_i) \right)^2}{ [\sigma^f(w_i)]^2 }  } \right)
\end{equation*}
where $n$ is the size of the validation set, $\hat{f}(w_i)$ is the prediction of our GP model at point $w_i$ and $f(w_i)$ is the associated true value and the validation set consists of 13000 evenly spaced points (see~\ref{subsec:cosine}). The values, reported in~\tabref{tab:resRoustant}, show that the PVA is constant, meaning that the estimation of the variance is kept proportional to the RMSE. The RMSE decreases as the number of hyperparameters is increasing. }

\begin{figure}[H]
\begin{center}
    \subfloat[EHH kernel: 79 hyperparameters, RMSE= 1.882.]{
    \centering
	\includegraphics[clip=true, height=4.5cm, width=7.2cm]{  HOMO_50_Roustant_curves.png}} \hspace{0.2 cm}
	\subfloat[FE kernel: 92 hyperparameters, RMSE= 22.610.]{
      \centering
	\includegraphics[clip=true, height=4.5cm, width=8cm]{  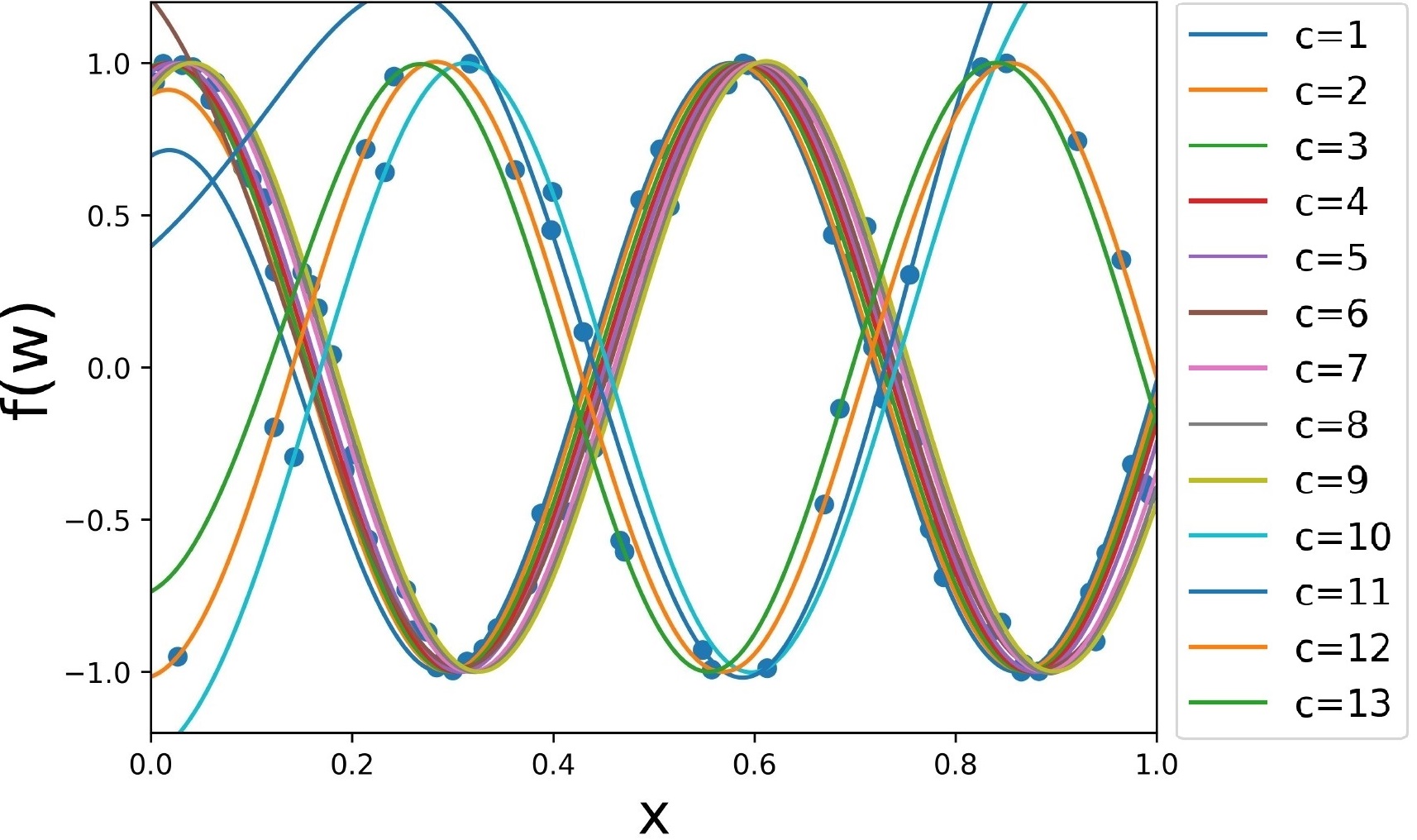}
	\label{fullroustant}} 	
\caption{Mean predictions comparison between EHH and FE kernels over the matrix $\Theta_1$ for the cosine problem with a 98 point DoE. }
\label{models_EHH}
\end{center} 
\end{figure}

\begin{table}[H]
\centering
 \caption{Kernel comparison for the cosine test case}
\begin{tabular}{ccccc}
  Kernel &  \# of Hyperparam. & RMSE & PVA  & CPU time (s) \\
  \hline 
  GD &  2 & 30.079&  21.99 & 1.4 \\   
  CR &  14 & 22.347 & 23.04& 24.5 \\
  EHH & 79 & 1.882 & 23.74 &514.5 \\
  HH  & 79 & 1.280   & 24.31 & 514.5 \\
\hline
\end{tabular}
\label{tab:resRoustant}
\end{table}
{\color{black}
\figref{models_EHH} shows a comparison between the FE and EHH kernels. Although the two kernels are equivalent in exact precision, the EHH kernel shows  more stable and better results in term of the accuracy compared to the FE general kernel. For this reason, in what comes next, only the EHH kernel will be considered on practical use cases.}

\begin{figure}[h]
\begin{center}

	\subfloat[GD kernel.]{
      \centering 
		\includegraphics[clip=true, height=4.5cm, width=5cm]{  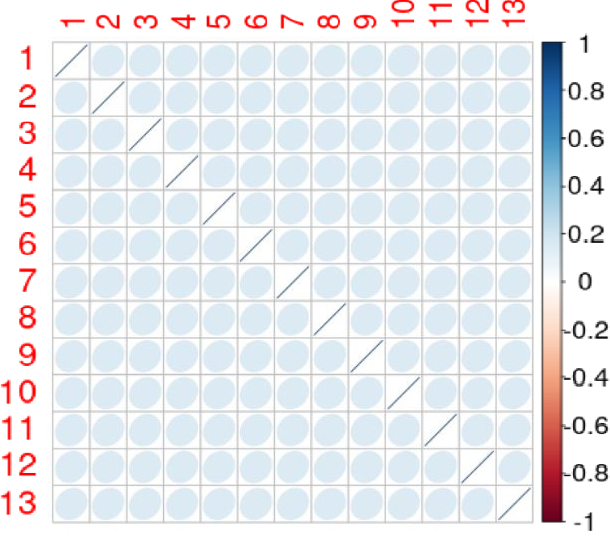} \label{corr_gower}
     }  
        \subfloat[CR kernel.]{
      \centering
		\includegraphics[clip=true, height=4.5cm, width=5cm]{  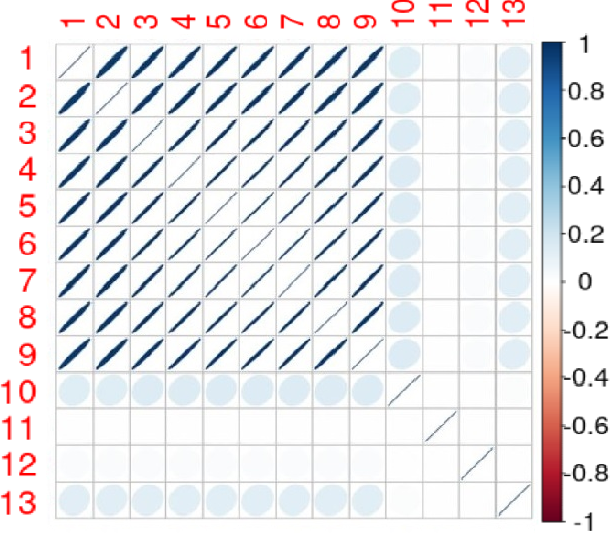} \label{corr_CR}
     }  
     
        \subfloat[EHH kernel.]{
      \centering
		\includegraphics[clip=true, height=4.5cm, width=5cm]{  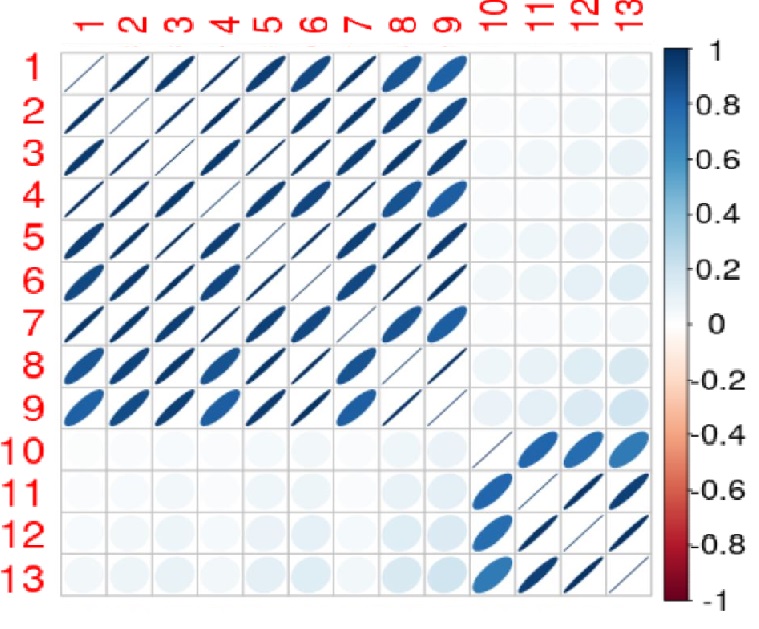}\label{corr_HS}
     }  
    \subfloat[HH kernel.]{
  \centering
	\includegraphics[clip=true, height=4.5cm, width=5cm]{  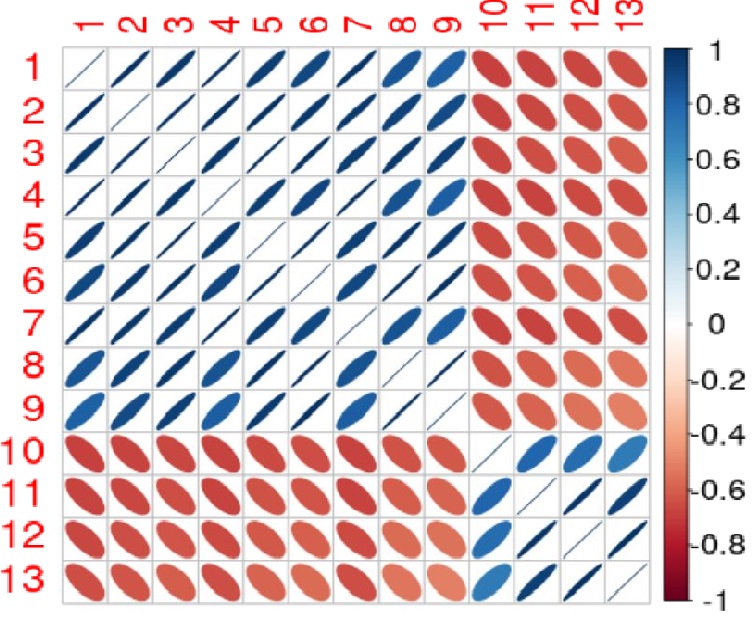}
    \label{corr_HH}
     }
\caption{Correlation matrix $R_1^{cat}$ using different choices for $\Theta_1$ on the cosine problem with a DoE of 98 points. }\label{corr_Roustant}
\end{center} 
\end{figure}

The estimated correlation matrix $R_i=R_1^{cat}$ is shown in~\figref{corr_Roustant}. For two given levels $\{\ell_r^1,\ell_s^1\}$, the correlation term $[R_1]_{\ell_r^1,\ell_s^1}$ is in blue for correlation values close to 1, in white for correlations close to 0 and in red for values close to -1; moreover the thinner the ellipse, the higher the correlation and we can see that the correlation between a level and itself is always 1. As expected, with GD kernel, there is only one estimated "mean correlation" as in~\figref{corr_gower}. For CR kernel (see~\figref{corr_CR}), the most important levels (1 to 9) are strongly correlated (in blue) with one another and the other levels (10 to 13) that should also have been correlated are badly estimated because of the kernel limitations that neglected them. 
In contrast, the EHH kernel (see~\figref{corr_HS}) gives a good approximation of the real correlations as it recovers the two groups of highly correlated levels. We recover the levels 1 to 9 as strongly similar and the levels 10 to 13 as strongly similar which is a good point but the two groups, even if less similar, are still positively correlated with one another. The latter between-group correlations should have been negative but the squared exponential kernel does not allow negative values. 
\textcolor{black}{
The comparison with the HH kernel, as proposed in~\cite{Zhou}, see~\figref{corr_HH}, shows that  even if the HH kernel is more general compared to EHH kernel, both kernels have an RMSE of the same order of magnitude (around 1.8 for EHH and 1.3 for HH).
}

\begin{figure}[H]
\centering
\subfloat[RMSE value versus DoE size.]{
\includegraphics[clip=true, height=5.3cm, width=7.65cm]{  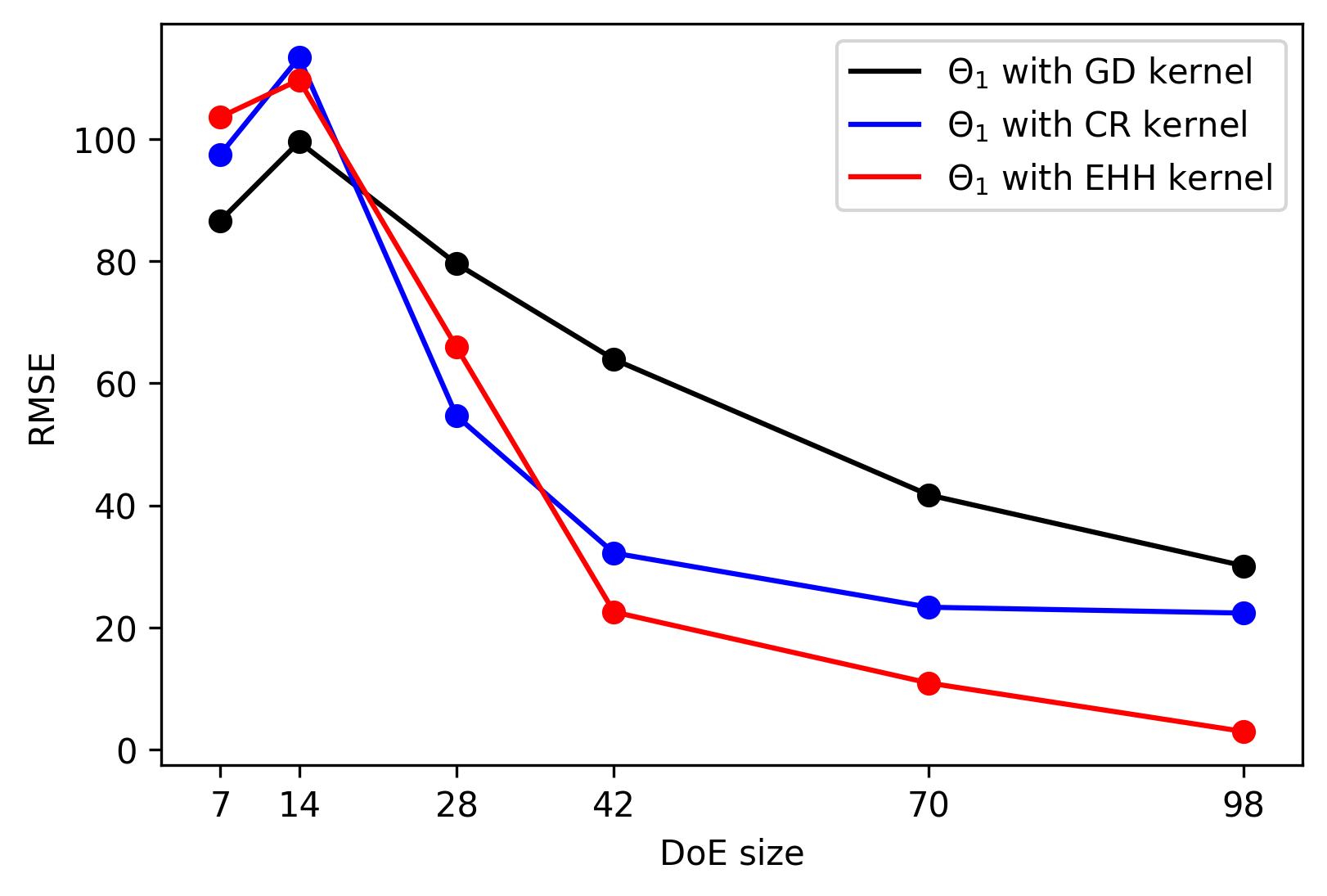}
}
\subfloat[CPU time (log scale) versus DoE size.]{
\includegraphics[clip=true, height=5.3cm, width=7.65cm]{  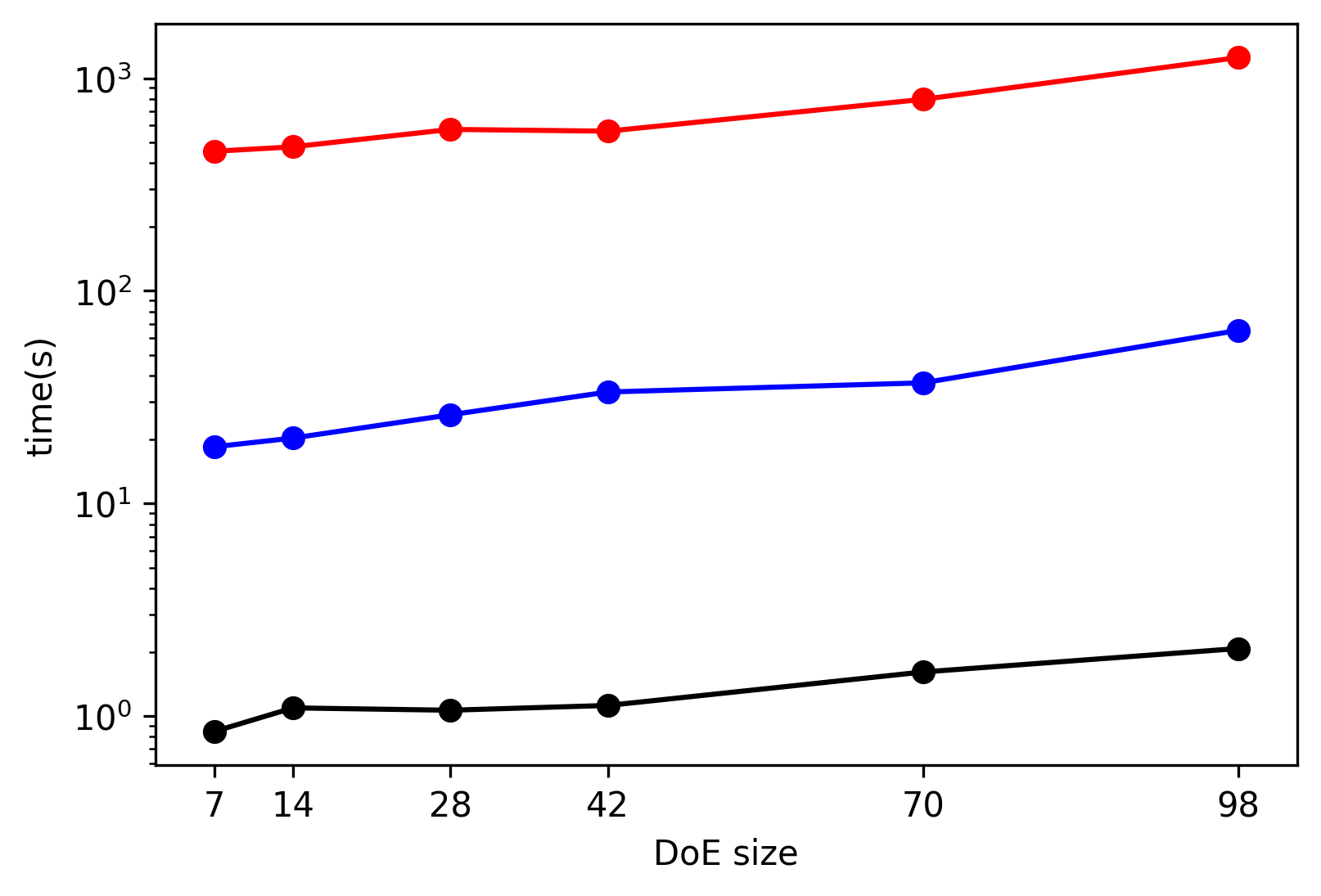}
\label{fig:time_increase}    
}
\caption{RMSE and CPU time to compute models with respect to DoE size.}
\label{fig:conv_RMSE}    
\end{figure}

On this particular test case, with a 98 point DoE, the more general the kernel, the better the performance and precision of the resulting GP. To show the DoE size impact, on~\figref{fig:conv_RMSE}, we draw 6 LHS DoEs of different sizes and we plot the RMSE and computational time for the kernels to see how they behave for different DoE sizes.

As expected, when the size of the DoE is too small for the problem (here smaller than 15 points), the three model behaviours are similarly bad because too little information is available for the hyperparameters optimization. However, when the size of the DoE is sufficiently large, we found the same hierarchy we found with 98 points on~\figref{models_Roustant} and the more complex the model, the faster the RMSE convergence. Nevertheless, on~\figref{fig:time_increase}, we can see that the computational costs of the models scale hardly with the DoE size on a logarithmic scale.

\subsection{Application to engineering problems}

To validate and compare our method on real applications, we will consider two engineering problems of different scale to analyze the model behaviour. In Section~\ref{subsec:beam}, we present an engineering beam bending problem and in Section~\ref{subsec:aircraft}, we introduce a complex system problem from aircraft design.

\subsubsection{Cantilever beam bending problem ($n=2$, $ m=0$, $l=1$ and $L_1=12$) }
\label{subsec:beam}
\begin{figure}[H]
\centering
\vspace{-10pt}
\subfloat[Bending problem]{
\begin{tikzpicture}
    \point{origin}{-0.75}{-0.25};

    \point{begin}{0}{0};
    \point{end}{5}{0};
    \point{end_bot}{5}{-1};
    \beam{2}{begin}{end};
    \support{3}{begin}[-90];
    \load{1}{end}[90]   ;
    \notation{1}{origin}{0};
    \notation{1}{end}{$F=50kN$};

     \draw[<->] (end) -- (end_bot) node[midway, right] {$\delta$} ;
     \draw[<->] (0,0.5) -- (5,0.5) node[midway, above] {L};
     
    \draw
      [-, ultra thick] (begin) .. controls (2.5, -.26) .. (5, -1)
      [-, ultra thick] (begin) .. controls (2.5, -.5) .. (5, -1.6)
      [-, ultra thick] (begin) .. controls (2.5, -.9)   .. (5, -2.2);
  \end{tikzpicture}
\label{fig:beam}   
}
\subfloat[Possible cross-section shapes.]{
\centering
\begin{tikzpicture}
\vspace{1cm}

\tstar{0.25}{0.5}{6}{0}{thick,fill=yellow}
\tstar{0.14}{0.28}{6}{0}{thick,fill=white}

\tstar{0.25}{0.5}{6}{0}{thick,fill=yellow,xshift=-1.2cm}
\tstar{0.08}{0.16}{6}{0}{thick,fill=white,xshift=-1.2cm}

\tstar{0.25}{0.5}{6}{0}{thick,fill=yellow,xshift=-2.4cm}

\def\pos{1.2}
\fill[black] (\pos-0.24,-0.5) -- (\pos-0.24,0.5) -- (\pos+0.24,0.5)  -- (\pos+0.24,-0.5)   -- cycle ;
\fill[black] (\pos-0.5,-0.5) -- (\pos-0.5,-0.18) -- (\pos+0.5,-0.18)  -- (\pos+0.5,-0.5)   -- cycle ; 
\fill[black] (\pos-0.5,0.18) -- (\pos-0.5,0.5) -- (\pos+0.5,0.5)  -- (\pos+0.5,0.18)   -- cycle ;  

\def\pos{2.4}
\fill[black] (\pos-0.19,-0.5) -- (\pos-0.19,0.5) -- (\pos+0.19,0.5)  -- (\pos+0.19,-0.5)   -- cycle ;
\fill[black] (\pos-0.5,-0.5) -- (\pos-0.5,-0.25) -- (\pos+0.5,-0.25)  -- (\pos+0.5,-0.5)   -- cycle ; 
\fill[black] (\pos-0.5,0.25) -- (\pos-0.5,0.5) -- (\pos+0.5,0.5)  -- (\pos+0.5,0.25)   -- cycle ;  

\def\pos{3.6}
\fill[black] (\pos-0.14,-0.5) -- (\pos-0.14,0.5) -- (\pos+0.14,0.5)  -- (\pos+0.14,-0.5)   -- cycle ;
\fill[black] (\pos-0.5,-0.5) -- (\pos-0.5,-0.32) -- (\pos+0.5,-0.32)  -- (\pos+0.5,-0.5)   -- cycle ; 
\fill[black] (\pos-0.5,0.32) -- (\pos-0.5,0.5) -- (\pos+0.5,0.5)  -- (\pos+0.5,0.32)   -- cycle ;  

\vspace{1.2cm}

\fill[green,even odd rule] (0,-1.2) circle (0.5) (0,-1.2) circle (0.33);
\draw (0,-1.2) circle (0.5) ;
\draw (0,-1.2) circle (0.33) 
; 
\fill[green,even odd rule] (-1.2,-1.2) circle (0.5)(-1.2,-1.2) circle (0.17);
\draw (-1.2,-1.2) circle (0.5) ;
\draw (-1.2,-1.2) circle (0.17) 
; 
\fill[green,even odd rule] (-2.4,-1.2) circle (0.5) ;
\draw (-2.4,-1.2) circle (0.5) ;

\def\pos{1.2}
\fill[blue,even odd rule]  (\pos-0.5,-1.7) -- (\pos-0.5,-0.7) -- (\pos+0.5,-0.7) -- (\pos+0.5,-1.7) -- cycle ;

\def\pos{2.4}
\fill[blue,even odd rule]  (\pos-0.5,-1.7) -- (\pos-0.5,-0.7) -- (\pos+0.5,-0.7) -- (\pos+0.5,-1.7) -- cycle   (\pos-0.25,-1.45) -- (\pos-0.25,-0.95) -- (\pos+0.25,-0.95) -- (\pos+0.25,-1.45) -- cycle ;

\def\pos{3.6}
\fill[blue,even odd rule]  (\pos-0.5,-1.7) -- (\pos-0.5,-0.7) -- (\pos+0.5,-0.7) -- (\pos+0.5,-1.7) -- cycle   (\pos-0.35,-1.55) -- (\pos-0.35,-0.85) -- (\pos+0.35,-0.85) -- (\pos+0.35,-1.55) -- cycle ;

  \end{tikzpicture}

\label{fig:beam_shape}    
}
\caption{Cantilever beam problem.}
\end{figure}
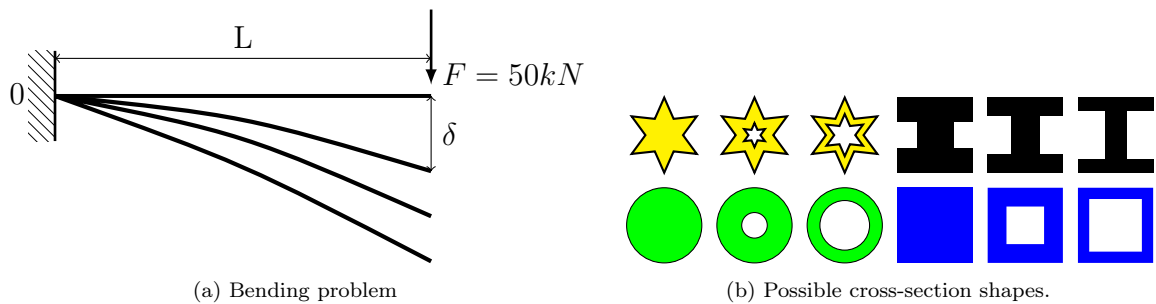

A first engineering problem commonly used for model validation is the beam bending problem in its linear elasticity range~\cite{Roustant, Cheng2015TrustRB}.
This problem is illustrated on~\figref{fig:beam} and consists of a cantilever beam loaded at its free extremity with a force $F$. As in Cheng $\textit{et al.}$~\cite{Cheng2015TrustRB}, we choose a constant Young modulus of $E=200$GPa and a load of $F=50$kN. Moreover, as in Roustant $\textit{et al.}$~\cite{Roustant}, we consider 12 possible cross-sections: there are 4 possible shapes, illustrated in~\figref{fig:beam_shape} that could be hollow, thick or full. For a given cross-section (shape and thickness), its size is determined by its surface $S$. Every cross section is associated with a normalized moment of inertia $\tilde{I}$ about the neutral axis. The latter is a latent variable associated to the beam shape~\cite{oune2021latent}.

Therefore, the problem to model has two continuous variables: the length $L \in [10,20]$ (in $m$) and the surface  $S \in [1,2]$ (in $m^2$) and one categorical variable $ \tilde{I}$ with 12 levels. The tip deflection, at the free end, $\delta$ is given by $$ \delta = f( \tilde{I}, L,S) = \frac{F}{3E} \frac{L^3}{S^2\tilde{I}}. $$



    
To compare our models, we draw a 98 point LHS as training set and the validation set is a grid of $12\times30\times30=10800$ points. For both squared exponential and absolute exponential kernels, the RMSE, likelihood and computational time for every model are shown in~\tabref{tab:resCantilever}. We recall that squared exponential and absolute exponential kernels differ only on the continuous variables and are the same for the categorical part.
As expected, the computational time and the likelihood increase when the model is more complex. The DoE seems of sufficient size for this problem as the computed RMSE (\textit{i.e.}, the total displacement error) decreases with the model complexity.

\begin{table}[H]
\centering
 \caption{Results of the cantilever beam models}
\small

\resizebox{0.9\columnwidth}{!}{%
\small

\begin{tabular}{ccccc}
  Categorical kernel & Continuous kernel & Displacement error (cm) & $ \ $ Likelihood  &$\ $ Time (s) \\
  \hline 
  GD & squared exponential &1.3858 & 111.13&  8.02 \\   
  CR & squared exponential & 1.1604 & 162.26 & 89.1 \\
  EHH & squared exponential &0.1247 &
256.90 & 2769.4 \\
  \hdashline 
  GD & absolute exponential & 3.2403 & 74.48   & 14.71  \\
  CR & absolute exponential & 3.0918 & 99.00 & 260.1 \\
  EHH & absolute exponential & 2.0951& 102.48 & 19784\\
\hline
\end{tabular}
}
\label{tab:resCantilever}
\end{table}

\begin{figure}[H]
\begin{center}
	\subfloat[With GD kernel.]{
      \centering 
\includegraphics[clip=true, height=4.6cm, width=5cm]{  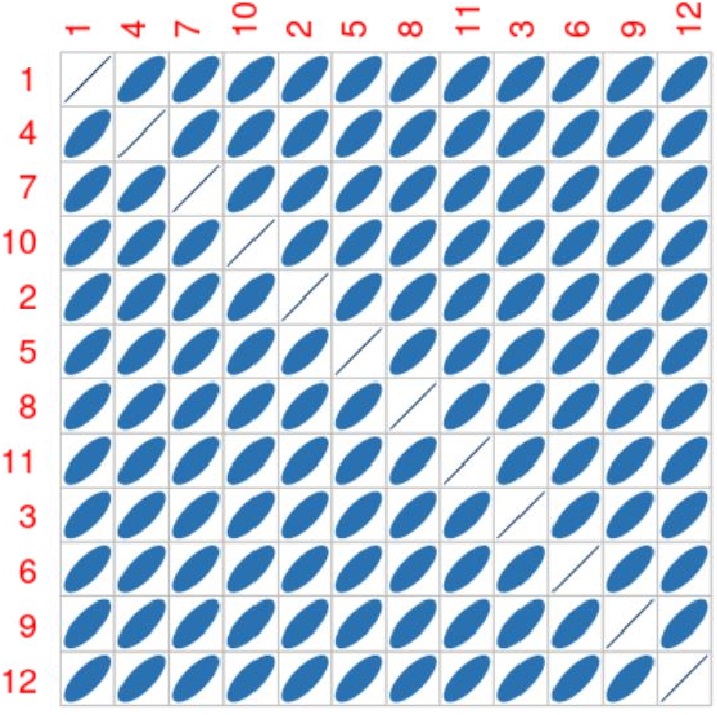} \label{corr_canti_gower}
     }  
	\subfloat[ With CR kernel.]{
      \centering
\includegraphics[clip=true, height=4.6cm, width=5cm]{  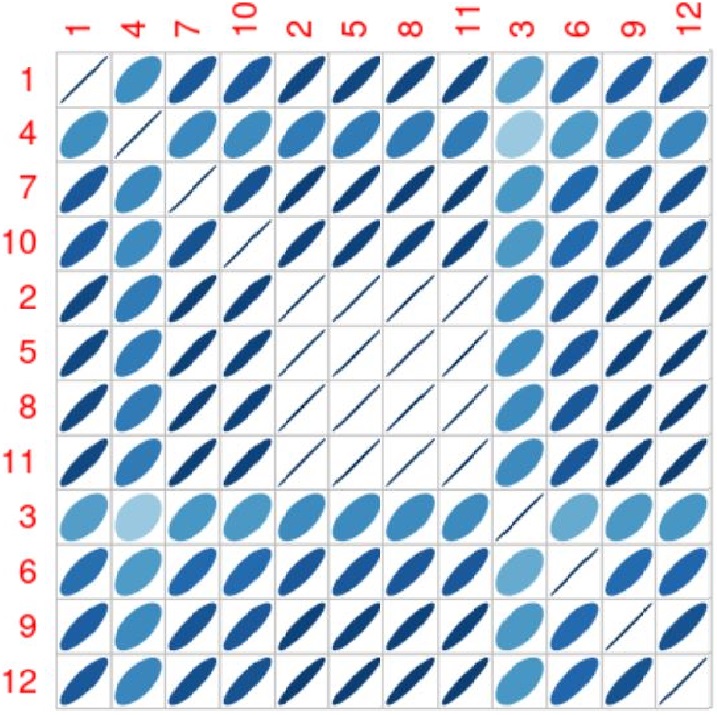} \label{corr_canti_cr}
     }  
	\subfloat[ With EHH kernel. ]{
      \centering 
\includegraphics[clip=true, height=4.6cm, width=5cm]{  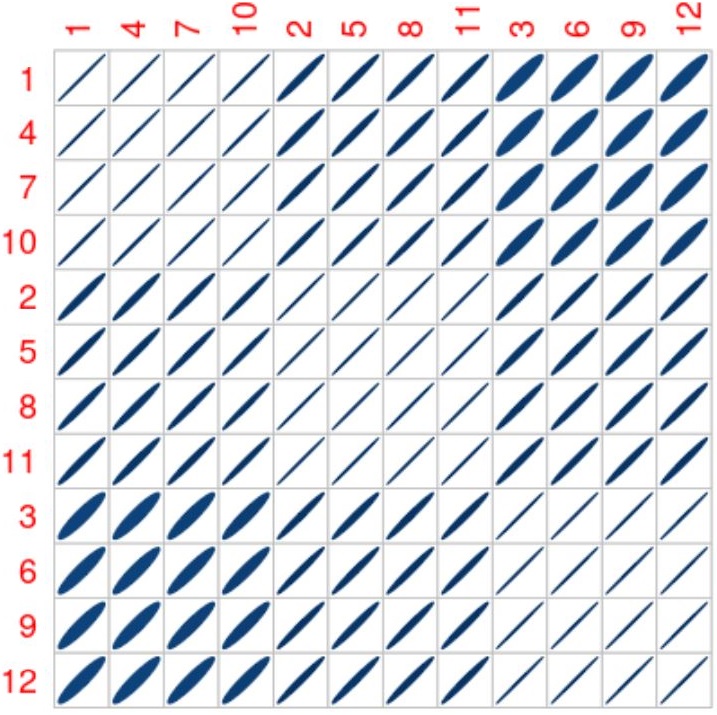} \label{corr_canti_ehh}
     }  
\includegraphics[clip=true, height=5cm, width=0.5cm]{  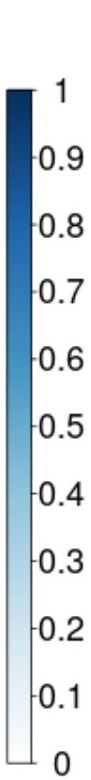}
\centering     
\caption{Correlation matrix $R_1^{cat}$  using different choices for $\Theta_1$ for the categorical variable $\tilde{I}$ from the cantilever beam problem.}
\label{corr_Cantilever}
\end{center} 
\end{figure}

In~\figref{corr_Cantilever}, we have drawn the correlation matrix found between the cross-section shape (the resulting $R_1$ correlation matrix) for the three models. On the figure below, the higher the correlation, the thinner the ellipse.

As expected, we have 3 groups of 4 shapes depending on their respective thickness (respectively, the levels \{1,4,7,10\} the levels \{2,5,8,11\} and the levels \{3,6,9,12\}). The more the thickness is similar, the higher the correlation: the thickness has more impact than the shape of the cross-section on the tip deflection. However, given the database, two points with similar $L$ and $S$ values will have similar output whatever the cross-section. The effect of the cross-section on the output is always the same (in the form of $\frac{1}{\tilde{I}}$) leading to an high correlation after maximizing the likelihood. In~\figref{corr_canti_ehh}, with the EHH kernel, we can distinguish the 3 groups of 4 shapes and, because the correlations are close to 1, the homoscedastic hyperphere model~\cite{Pelamatti} would lead to the same correlation matrix.  Also, with the CR kernel of~\figref{corr_canti_cr}, the medium thick group \{2,5,8,11\} being correlated with both the full and the hollow group, its correlation values are the higher whereas the correlation hyperparameters associated to the two other groups are smaller. 
For the GD model in~\figref{corr_canti_gower}, there is only one mean positive correlation value as before.

\subsubsection{Aircraft design application ($n=10$, $ m=0$, $l=2$ and $L_1=9$, $L_2=2$) }
\label{subsec:aircraft}

The ``\texttt{DRAGON}'' aircraft concept 
has been introduced by ONERA in 2019~\cite{schmollgruber} within the scope of the European CleanSky 2 program\footnote{\href{https://www.cleansky.eu/technology-evaluator}{\color{blue}https://www.cleansky.eu/technology-evaluator}} which sets the objective of 30\% reduction of CO2 emissions by 2035 with respect to 2014 state of the art.
The employment of a distributed propulsion comes at a certain cost; a turboelectric propulsive chain is necessary to power the electric fans which brings additional complexity and weight.
The turboelectric propulsive chain being an important weight penalty, it is of particular interest to optimize the chain and particularly the number and type of each component, characterized by some discrete values.  The definition of the architecture variable is given in~\tabref{tab:dragon_archi1} and the definition of the turboshaft layout is given in~\tabref{tab:dragon_archi2}. For the sake of simplicity, we restrict the optimization problem to the case of two electric cores and generators but more optimizations have been performed in~\cite{SciTech_cat}.

%
%

\begin{table}[H]
\centering
\vspace*{-0.1cm}

 \subfloat[Definition of the architecture variable and its 9 associated levels.]{
\small

\resizebox{0.8\columnwidth}{!}{%
\small

\begin{tabular}{cccc}
  Architecture number $\ $ & Number of motors $\ $ & Number of cores $\ $ & Number of generators $\ $ \\
  \hline
  1 & 8 &2 & 2\\
  2 & 12 & 2 & 2\\
  3 & 16 & 2 & 2\\
  4 &20 &2 & 2\\
  5 & 24 & 2 & 2\\
  6 & 28 & 2 & 2\\
  7 &32 & 2 & 2\\
  8 & 36  & 2 & 2\\
  9 & 40 & 2 & 2\\

\hline
\end{tabular}
}
\label{tab:dragon_archi1}
}
\centering
\vspace*{+0.1cm}

\subfloat[Definition of the turboshaft layout variable and its 2 associated levels.]{
\small

\resizebox{0.8\columnwidth}{!}{%
\small

\begin{tabular}{cccccc}
  Layout & Position & y ratio & Tail & VT aspect ratio & VT taper ratio\\
  \hline 
  1 & under wing &0.25 & without T-tail& 1.8 & 0.3 \\
  2 & behind & 0.34 & with T-tail& 1.2 & 0.85\\
 
\hline
\end{tabular}
}
\label{tab:dragon_archi2}
}
\caption{Categorical variable definition}
\end{table}
The analysis of ``\texttt{DRAGON}'' is treated with Overall Aircraft Design method in FAST-OAD~\cite{David_2021}.  We are considering the following problem described in~\tabref{tab:dragon}.
\begin{table}[H]
\centering
\vspace*{-0.3cm}

 \caption{Definition of the ``\texttt{DRAGON}'' optimization problem.}
\small

\resizebox{1.0\columnwidth}{!}{%
\small

\begin{tabular}{lllrr}
 & Function/variable & Nature & Quantity & Range\\
\hline
\hline
Model & Fuel mass & cont & 1 &\\
\hline
with respect to & \mbox{Fan operating pressure ratio} & cont & 1 & $\left[1.05, 1.3\right]$ \\  
     & \mbox{Wing aspect ratio} & cont & 1 &    $\left[8, 12\right]$ \\
    & \mbox{Angle for swept wing} & cont & 1 & $\left[15, 40\right]$  ($^\circ$) \\
     & \mbox{Wing taper ratio} & cont & 1 &    $\left[0.2, 0.5\right]$ \\
     & \mbox{HT aspect ratio} & cont & 1 &    $\left[3, 6\right]$ \\
    & \mbox{Angle for swept HT} & cont & 1 & $\left[20, 40\right]$  ($^\circ$) \\
     & \mbox{HT taper ratio} & cont & 1 &    $\left[0.3, 0.5\right]$ \\
 & \mbox{TOFL for sizing}  & cont &1 & $\left[1800, 2500\right]$ ($m$) \\
 & \mbox{Top of climb vertical speed for sizing $ \ $} & cont & 1 & $\left[300, 800\right]$ ($ft/min$) \\
 & \mbox{Start of climb slope angle} & cont & 1 & $\left[0.075, 0.15\right] $ ($rad$) \\

 & \multicolumn{2}{l}{Total  continuous variables} & 10 & \\
 \cline{2-5}
& \mbox{Architecture} & cat & 9 levels & \{1,2,3, \ldots,7,8,9\} \\
& \mbox{Turboshaft layout} & cat & 2 levels & \{1,2\} \\

 & \multicolumn{2}{l}{Total categorical variables} & 2 & \\
 \cline{2-5}

  &   \multicolumn{2}{l}{\textbf{Total relaxed variables}} & {\textbf{21}} & \\
  \hline

\end{tabular}
}
\label{tab:dragon}
\end{table}

\begin{figure}[H]
\begin{center}
	\subfloat[GD kernel.]{
      \centering 
\includegraphics[clip=true, height=4.4cm, width=5cm]{  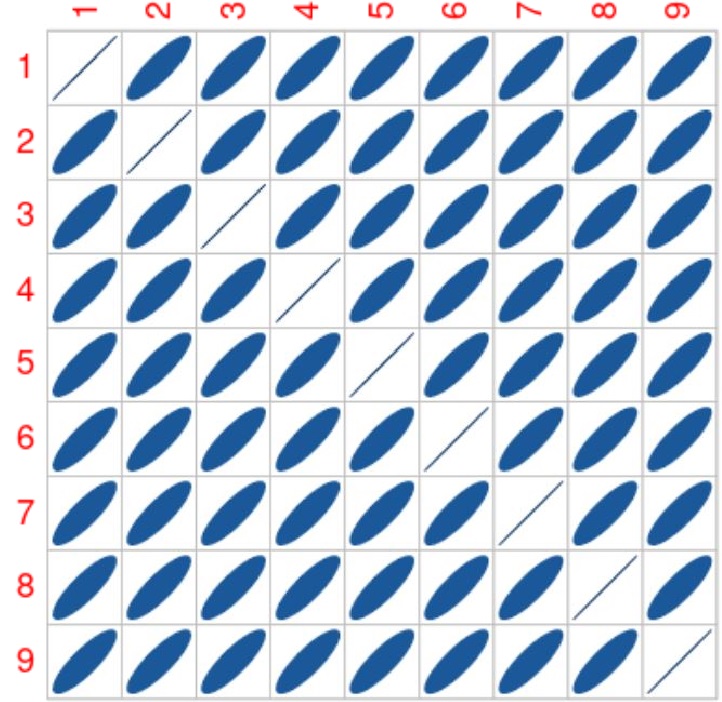} \label{corr_turboelectric_gower}
     }  
	\subfloat[CR kernel.]{
      \centering 
\includegraphics[clip=true, height=4.4cm, width=5cm]{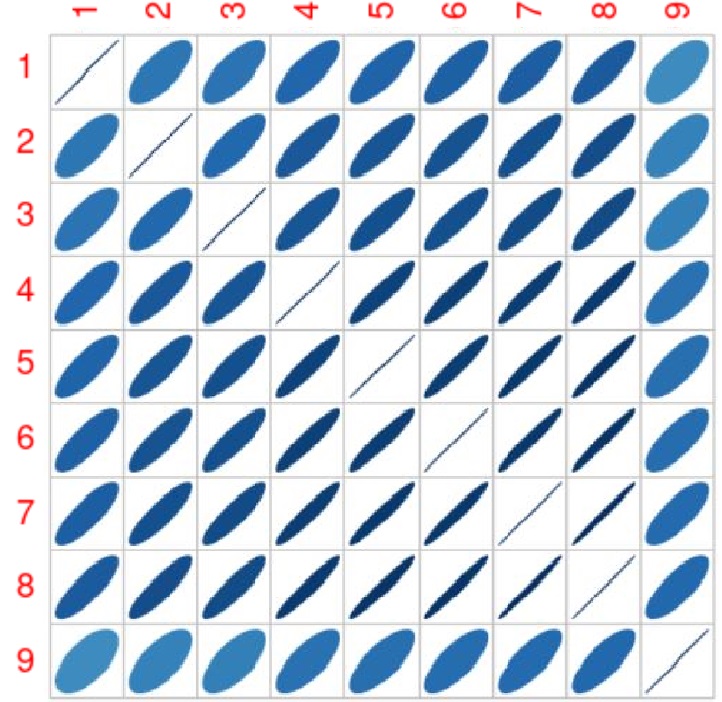} \label{corr_turboelectric_cr}
     }  
	\subfloat[EHH kernel.]{
      \centering 
\includegraphics[clip=true,  height=4.4cm, width=5cm]{  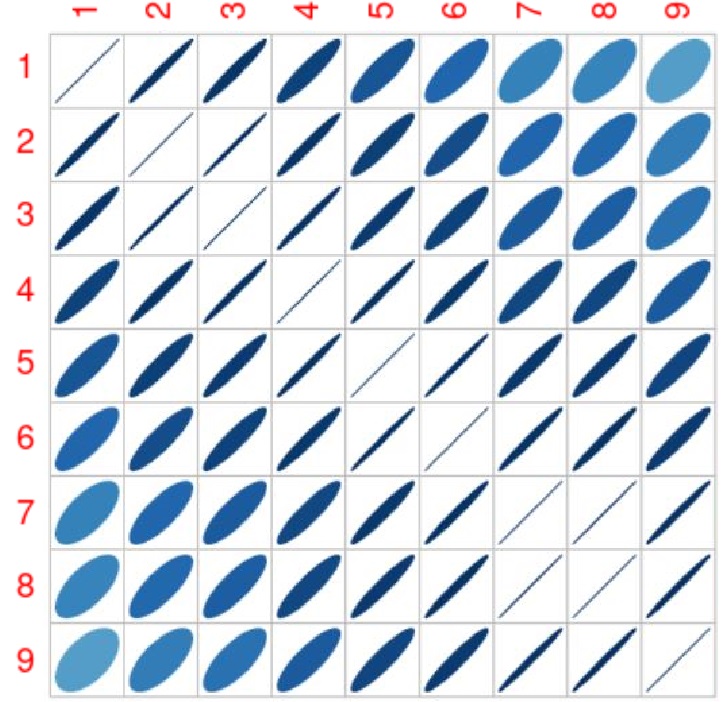}
\label{corr_turboelectric_ehh}
     }  
\includegraphics[clip=true, height=4.7cm, width=0.5cm]{  leg.jpg}
\centering     
\caption{Correlation matrix $R_1^{cat}$  using different choices for $\Theta_1$  for the turboelectric architecture variable.}
\label{corr_turboelectric}
\end{center} 
\end{figure}

\begin{figure}[H]
\begin{center}
\vspace{-1cm}
	\subfloat[GD kernel.]{
      \centering 
\includegraphics[clip=true, height=4.4cm, width=5cm]{  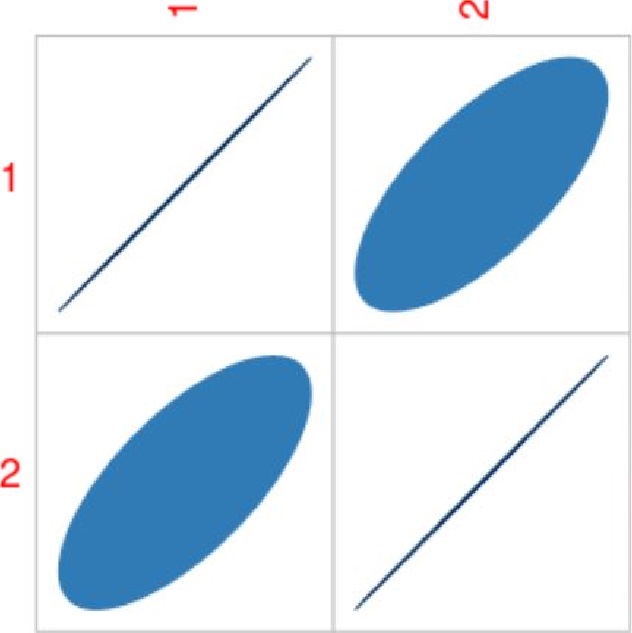} \label{corr_turboshaft_gower}
     }  
	\subfloat[CR kernel.]{
      \centering 
\includegraphics[clip=true, height=4.4cm, width=5cm]{  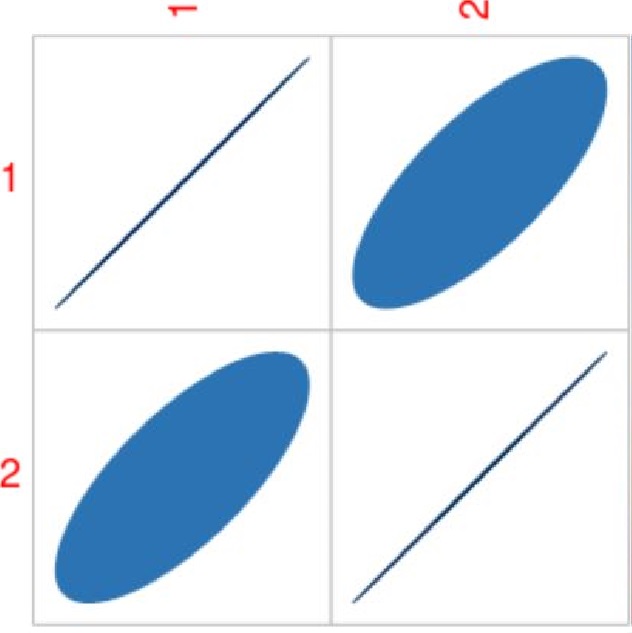} \label{corr_turboshaft_cr}
     }  
	\subfloat[EHH kernel.]{
      \centering 
\includegraphics[clip=true, height=4.4cm, width=5cm]{  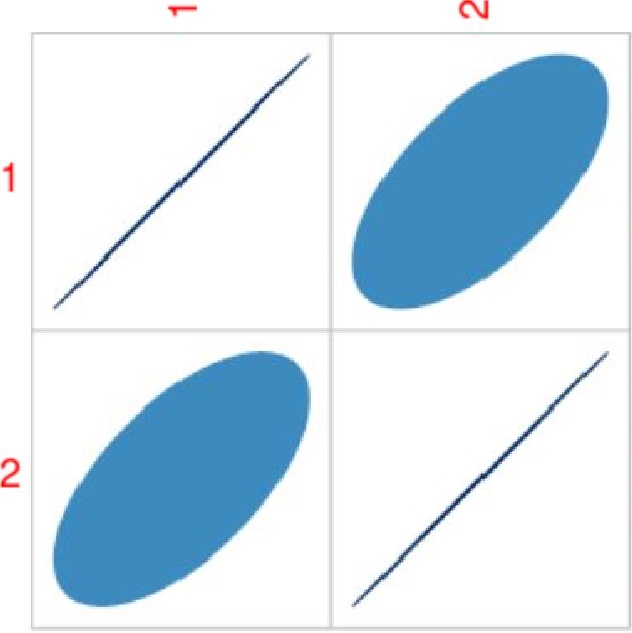} \label{corr_turboshaft_ehh}
     }  
\centering
\includegraphics[clip=true, height=4.65cm, width=0.5cm]{  leg.jpg}
\caption{Correlation matrix $R_2^{cat}$  using different choices for $\Theta_2$  for the turboshaft layout variable.}
\label{corr_turboshaft}
\end{center} 
\end{figure}

Twice, we draw  $250$ points by LHS. Over the first DoE, that is the training set, we build the model to predict the fuel mass and over the second one, we validate our prediction and compute the RMSE reported in~\tabref{tab:resDragon}.
In this case, the number of hyperparameters is 12 for GD kernel, 21 for CR kernel and 47 for EHH kernel. Evaluating the function is costly, around 4 minutes for a single point. We observed similar performances for all models, the performance is mostly determined by the choice of the continuous kernel. 
For a problem that has that many variables, it seems useless and impractical to use a complicated model, the GD kernel being already performing well. 
On~\figref{corr_turboelectric}, we plot, for the three kernels, the approximate correlation matrices for the first categorical variable. As we can see, when considering the general EHH kernel, as in~\figref{corr_turboelectric_ehh}, the closer the levels, the higher the correlation. In fact, in this case, the only difference between two levels is the number of motors. Therefore, the more similar the number of motors, the more similar the fuel consumption. Given that, we expect, when considering CR kernel as in~\figref{corr_turboelectric_cr} that the higher correlation should appear "in the middle" \{4,5,6\} as these levels are meant to be the most correlated with the others. This is what happens to a certain extent but the levels 7 and 8 are weirdly appearing too much correlated with one another. This could be a numerical problem, the optimization being hard with that many variables and hyperparameters. As before, the GD kernel is the less precise and just give a mean correlation over the whole space as in~\figref{corr_turboelectric_gower}. In~\figref{corr_turboshaft}, we plot, for the three methods, the approximated correlation matrices for the second categorical variable. There is only two engine layouts so there is only one correlation. In this case, the correlation is positive indicating that the plane behave in the same way no matter the layout.

\begin{table}[H]    
\centering
 \caption{Results of the aircraft models based on a 250 point validation set}
\small
\resizebox{0.8\columnwidth}{!}{%
\small

\begin{tabular}{ccccc}
  kernel& $ \ $ number of hyperparameters  & $\ $  kernel   & $\ $fuel error (kg) $\ $  & time (s)   \\
  \hline 
  GD & 12 & squared exponential & 2115 & 65 \\
  CR & 21 & squared exponential & 2068  & 210  \\
  EHH & 47 & squared exponential & 2147 &  9450\\
   \hdashline 
  GD & 12 & absolute exponential &1666 & 65 \\
  CR & 21 & absolute exponential & 1664  & 210 \\
  EHH & 47 & absolute exponential & 1593 &  9295 \\
\hline
\end{tabular}
}
\label{tab:resDragon}
\end{table}

One can note that increasing the number of motors or changing a layout will not change the way an aircraft flies. For example, having more motors will only increase the fuel consumption by a given factor. The latter will always remain positive and related to the continuous variables. Hence, in this test case, we do not have opposite effects between two categorical levels.

In most industrial applications, radically opposite effects over a complex system do not occur so often. For instance, on the industrial applications that can be found on the literature, there was not a clear need for negative correlation values~\cite{Pelamatti,Roustant, cuesta2021comparison}. Therefore, in practice, the exponential model is not that limiting compared to the homoscedastic hypersphere model.

\section{Conclusion}
\label{sec:conclu}

In this work, we have proposed a class of kernels for GP models that extends the exponential continuous kernels to the mixed-categorical setting. We showed that this class of kernels generalizes Gower distance and continuous relaxation based kernels. A classification between the proposed kernels as well as a proof of the SPD nature of the resulting correlation matrices have been also proposed. Numerical illustrations on analytical toy problems showed the good potential of the proposed kernels to reduce the number of hyperparameters and thus the computational time. The implementation of our proposed method has been released in the toolbox SMT v2.0\footnote{\url{https://smt.readthedocs.io/en/latest/}}. 

When considering complex kernels, a good approach would be to use a model reduction technique such as Kriging with Partial Least Squares (KPLS)~\cite{Bouhlel18} that is derived from the construction of the correlation matrix via a kernel function. KPLS is an adaptation of the Partial Least Squares regression for exponential kernels and is used to reduce the number of hyperparameters and handle a large number of mixed inputs. Further works will consider to include such dimension reduction techniques to improve the computational efficiency of our model  and tackle higher dimensional problems.

\section*{Acknowledgements}
This work is part of the activities of ONERA - ISAE - ENAC joint research group. The research presented in this paper has been performed in the framework of the AGILE 4.0 project (Towards Cyber-physical Collaborative Aircraft Development) and has received funding from the European Union Horizon 2020 Programme under grant agreement n${^\circ}$ 815122. 
We thank Raul Carreira Rufato (ISAE-Supaero MSc) for his contribution to Gower distance implementation, and Dr. Eric Nguyen Van (ONERA) and Christophe David (ONERA) for their contribution to DRAGON aircraft design.
\textcolor{black}{We also thanks the collaborators of SMT, namely ONERA, ISAE-SUPAERO, ICA (CNRS), NASA Glenn, the University of Michigan, Polytechnique Montréal and the University of California San Diego and in particular Rémi Lafage (ONERA).}
The authors are grateful to the partners of the AGILE 4.0 consortium for their contribution and feedback.




\appendix
\addcontentsline{toc}{section}{Appendix}
\section*{Appendix}

In~\ref{apendix:EHH2CR}, we give the parameterization that allows us to obtain the continuous relaxation kernel using our proposed framework. 
In~\ref{subsec:cosine}, the cosine test case is detailed. 

\section{Continuous relaxation is a particular instance of our proposed FE Kernel.} 
\label{apendix:EHH2CR}
To show that CR is a particular instance of FE, it suffices to show that the matrix $\Phi(\Theta_i)$ is diagonal whenever $\Theta_i$ is set to a diagonal one. In fact, assume that we have, in our general model, $[\Theta_i]_{j \neq j'} = 0, \  \forall (j,j') \in \{ 1, \ldots, L_i \}.$ 
Knowing that $\cos(0)=1$ and $\sin(0)=0$, the matrix $C(\Theta_i)$ writes as  \\
$$C(\Theta_i) = \begin{bmatrix}
1 & 0 & 0  & 0  \\
1  & 0 &  \ldots & 0 \\
\vdots &\vdots & \ddots & 0 \\
1 &  0 & 0  & 0 \\
\end{bmatrix} \quad \mbox{ and } \quad  C(\Theta_i) C(\Theta_i)^\top = \begin{bmatrix}
1 & 1 & 1  & 1  \\
1  & 1 &  \ldots & 1 \\
\vdots &\vdots & \ddots & 1 \\
1 &  1 & 1  & 1 \\
\end{bmatrix} $$
Therefore, we also have 
$$[\Phi(\Theta_i)]_{j \neq j'} = \frac{\log \epsilon }{2} ([C(\Theta_i) C(\Theta_i)^\top]_{j,j'} -1) = 0, \  \forall (j,j') \in \{ 1, \ldots, L_i \} $$ that is the continuous relaxation kernel. \qed \mbox{}\\ 


\section{Categorical cosine case}
\label{subsec:cosine}
\begin{figure}[H]
\centering
\includegraphics[scale=.12]{  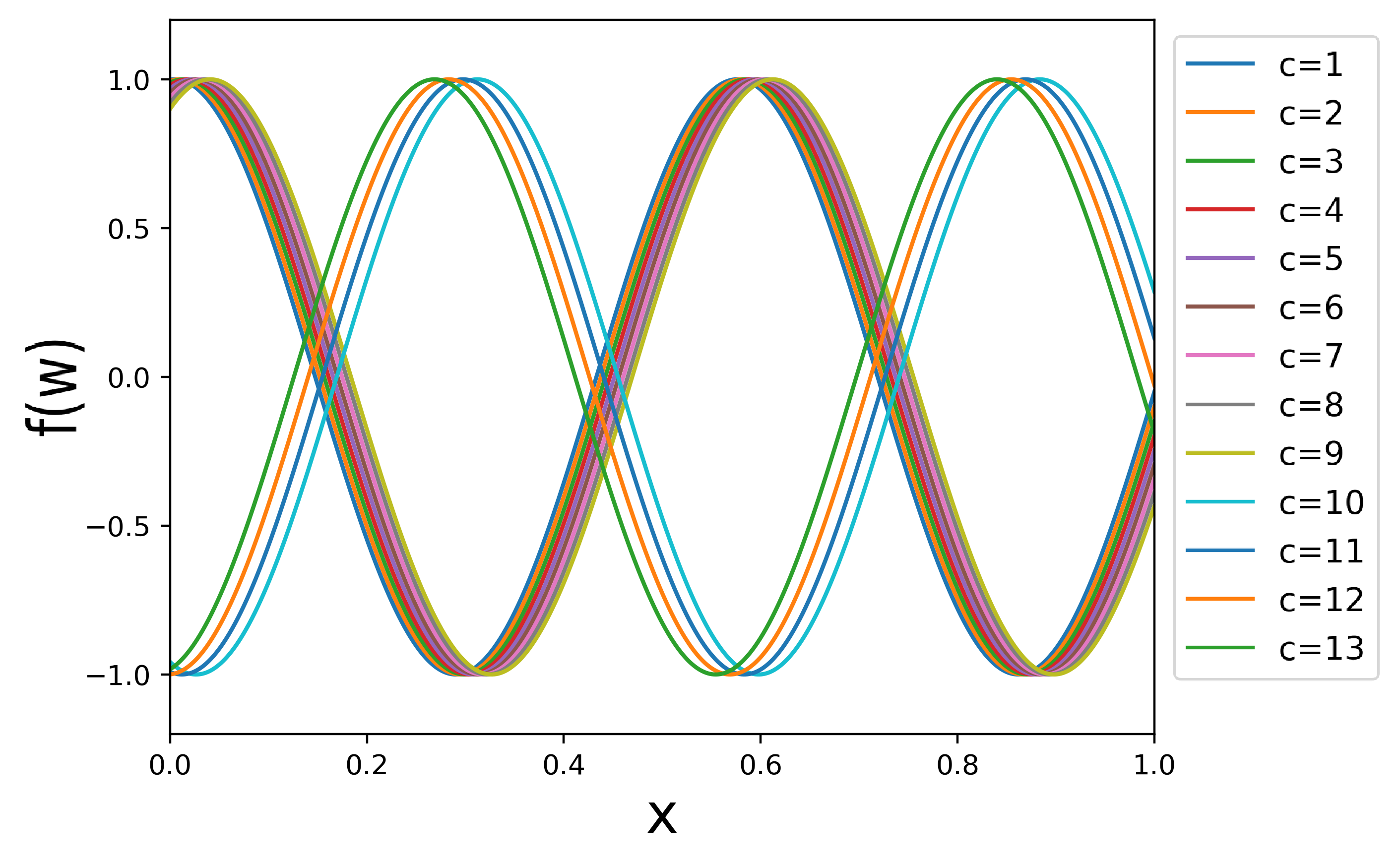}
\caption{Landscape of the cosine test case from~\cite{Roustant}.}
\label{fig:Roustant_ref}    
\end{figure}

This test case has one categorical variable with 13 levels and one continuous variable in $[0,1]$~\cite{Roustant}.
Let $w= (x,c )$ be a given point with  $x$ being the continuous variable and $c$ being the categorical variable, $c \in \{1, \ldots, 13\}$.

\begin{equation*}
\begin{split}
f(w) &= \cos \left( \frac{7 \pi}{2} x + \left( 0.4 \pi  + \frac{\pi }{15} c  \right) - \frac{c}{20} \right) , ~~~\mbox{if c $\in\{10,\ldots,9\}$ }  \\
f(w) &= \cos \left( \frac{7 \pi}{2} x  - \frac{c}{20} \right) , ~~~\mbox{if c $\in\{10,\ldots,13\}$ }  \\
\end{split}
\end{equation*}
The reference landscapes of the objective function (with respect to the categorical choices) are drawn on~\figref{fig:Roustant_ref}.

The DoE is given by a LHS of 98 points.
Our validation set is a evenly spaced grid of 1000 points in $x$ ranging  for every of the 13 categorical levels  for a total of 13000 points.


\bibliographystyle{elsarticle-num-names} 
\bibliography{main}

\begin{thebibliography}{57}
\expandafter\ifx\csname natexlab\endcsname\relax\def\natexlab#1{#1}\fi
\providecommand{\url}[1]{\texttt{#1}}
\providecommand{\href}[2]{#2}
\providecommand{\path}[1]{#1}
\providecommand{\DOIprefix}{doi:}
\providecommand{\ArXivprefix}{arXiv:}
\providecommand{\URLprefix}{URL: }
\providecommand{\Pubmedprefix}{pmid:}
\providecommand{\doi}[1]{\href{http://dx.doi.org/#1}{\path{#1}}}
\providecommand{\Pubmed}[1]{\href{pmid:#1}{\path{#1}}}
\providecommand{\bibinfo}[2]{#2}
\ifx\xfnm\relax \def\xfnm[#1]{\unskip,\space#1}\fi
\bibitem[{Saves et~al.(2022)Saves, {Nguyen Van}, Bartoli, Diouane, Lefebvre,
  David, Defoort, and Morlier}]{SciTech_cat}
\bibinfo{author}{P.~Saves}, \bibinfo{author}{E.~{Nguyen Van}},
  \bibinfo{author}{N.~Bartoli}, \bibinfo{author}{Y.~Diouane},
  \bibinfo{author}{T.~Lefebvre}, \bibinfo{author}{C.~David},
  \bibinfo{author}{S.~Defoort}, \bibinfo{author}{J.~Morlier},
\newblock \bibinfo{title}{Bayesian optimization for mixed variables using an
  adaptive dimension reduction process: applications to aircraft design},
\newblock in: \bibinfo{booktitle}{AIAA SciTech 2022 Forum},
  \bibinfo{year}{2022}.
\bibitem[{Snoek et~al.(2015)Snoek, Rippel, Swersky, Kiros, Satish, Sundaram,
  Patwary, Prabhat, and Adams}]{snoek2015scalable}
\bibinfo{author}{J.~Snoek}, \bibinfo{author}{O.~Rippel},
  \bibinfo{author}{K.~Swersky}, \bibinfo{author}{R.~Kiros},
  \bibinfo{author}{N.~Satish}, \bibinfo{author}{N.~Sundaram},
  \bibinfo{author}{M.~Patwary}, \bibinfo{author}{M.~Prabhat},
  \bibinfo{author}{R.~Adams},
\newblock \bibinfo{title}{Scalable bayesian optimization using deep neural
  networks},
\newblock in: \bibinfo{booktitle}{International conference on machine
  learning}, \bibinfo{year}{2015}.
\bibitem[{López-Lopera et~al.(2022)López-Lopera, Idier, Rohmer, and
  Bachoc}]{lopez}
\bibinfo{author}{A.~F. López-Lopera}, \bibinfo{author}{D.~Idier},
  \bibinfo{author}{J.~Rohmer}, \bibinfo{author}{F.~Bachoc},
\newblock \bibinfo{title}{Multioutput gaussian processes with functional data:
  A study on coastal flood hazard assessment},
\newblock \bibinfo{journal}{Reliability Engineering \& System Safety}
  \bibinfo{volume}{218} (\bibinfo{year}{2022}) \bibinfo{pages}{108139}.
\bibitem[{Ghasemi et~al.(2021)Ghasemi, Karbasi, {Zamani Nouri}, {Sarai
  Tabrizi}, and Azamathulla}]{MLP}
\bibinfo{author}{P.~Ghasemi}, \bibinfo{author}{M.~Karbasi},
  \bibinfo{author}{A.~{Zamani Nouri}}, \bibinfo{author}{M.~{Sarai Tabrizi}},
  \bibinfo{author}{H.~Z. Azamathulla},
\newblock \bibinfo{title}{Application of gaussian process regression to
  forecast multi-step ahead spei drought index},
\newblock \bibinfo{journal}{Alexandria Engineering Journal}
  \bibinfo{volume}{60} (\bibinfo{year}{2021}) \bibinfo{pages}{5375--5392}.
\bibitem[{Perera et~al.(2014)Perera, Anderson, and Ghosh}]{retina}
\bibinfo{author}{N.~Perera}, \bibinfo{author}{R.~Anderson},
  \bibinfo{author}{B.~Ghosh},
\newblock \bibinfo{title}{Detection of moving targets in the visual pathways of
  turtles using computational models},
\newblock in: \bibinfo{booktitle}{7th International Conference on Information
  and Automation for Sustainability}, \bibinfo{year}{2014}.
\bibitem[{Diouane et~al.(2016)Diouane, Gratton, Vasseur, Vicente, and
  Calandra}]{YDiouane_SGratton_XVasseur_LNVicente_HCalandra_2016}
\bibinfo{author}{Y.~Diouane}, \bibinfo{author}{S.~Gratton},
  \bibinfo{author}{X.~Vasseur}, \bibinfo{author}{L.~N. Vicente},
  \bibinfo{author}{H.~Calandra},
\newblock \bibinfo{title}{A parallel evolution strategy for an earth imaging
  problem in geophysics},
\newblock \bibinfo{journal}{Optim. Eng.} \bibinfo{volume}{17}
  (\bibinfo{year}{2016}) \bibinfo{pages}{3--26}.
\bibitem[{Rufato et~al.(2022)Rufato, Diouane, Henry, Ahlfeld, and
  Morlier}]{RaulAIAA}
\bibinfo{author}{R.~C. Rufato}, \bibinfo{author}{Y.~Diouane},
  \bibinfo{author}{J.~Henry}, \bibinfo{author}{R.~Ahlfeld},
  \bibinfo{author}{J.~Morlier},
\newblock \bibinfo{title}{A mixed-categorical data-driven approach for
  prediction and optimization of hybrid discontinuous composites performance},
\newblock in: \bibinfo{booktitle}{AIAA AVIATION 2022 Forum},
  \bibinfo{year}{2022}.
\bibitem[{Williams and Rasmussen(2006)}]{williams2006gaussian}
\bibinfo{author}{C.~K. Williams}, \bibinfo{author}{C.~E. Rasmussen},
  \bibinfo{title}{Gaussian processes for machine learning},
  \bibinfo{publisher}{MIT press Cambridge, MA}, \bibinfo{year}{2006}.
\bibitem[{Krige(1951)}]{krige1951statistical}
\bibinfo{author}{D.~G. Krige},
\newblock \bibinfo{title}{A statistical approach to some basic mine valuation
  problems on the witwatersrand},
\newblock \bibinfo{journal}{Journal of the Southern African Institute of Mining
  and Metallurgy} \bibinfo{volume}{52} (\bibinfo{year}{1951})
  \bibinfo{pages}{119--139}.
\bibitem[{Pelamatti et~al.(2019)Pelamatti, Brevault, Balesdent, Talbi, and
  Guerin}]{Pelamatti}
\bibinfo{author}{J.~Pelamatti}, \bibinfo{author}{L.~Brevault},
  \bibinfo{author}{M.~Balesdent}, \bibinfo{author}{E.-G. Talbi},
  \bibinfo{author}{Y.~Guerin},
\newblock \bibinfo{title}{Efficient global optimization of constrained mixed
  variable problems},
\newblock \bibinfo{journal}{Journal of Global Optimization}
  \bibinfo{volume}{73} (\bibinfo{year}{2019}) \bibinfo{pages}{583--613}.
\bibitem[{Zhou et~al.(2011)Zhou, Qian, and Zhou}]{Zhou}
\bibinfo{author}{Q.~Zhou}, \bibinfo{author}{P.~Z.~G. Qian},
  \bibinfo{author}{S.~Zhou},
\newblock \bibinfo{title}{A simple approach to emulation for computer models
  with qualitative and quantitative factors},
\newblock \bibinfo{journal}{Technometrics} \bibinfo{volume}{53}
  (\bibinfo{year}{2011}) \bibinfo{pages}{266--273}.
\bibitem[{Deng et~al.(2017)Deng, Lin, Liu, and Rowe}]{Deng}
\bibinfo{author}{X.~Deng}, \bibinfo{author}{C.~D. Lin},
  \bibinfo{author}{K.~Liu}, \bibinfo{author}{R.~K. Rowe},
\newblock \bibinfo{title}{Additive gaussian process for computer models with
  qualitative and quantitative factors},
\newblock \bibinfo{journal}{Technometrics} \bibinfo{volume}{59}
  (\bibinfo{year}{2017}) \bibinfo{pages}{283--292}.
\bibitem[{Roustant et~al.(2020)Roustant, Padonou, Deville, Cl{é}ment, Perrin,
  Giorla, and Wynn}]{Roustant}
\bibinfo{author}{O.~Roustant}, \bibinfo{author}{E.~Padonou},
  \bibinfo{author}{Y.~Deville}, \bibinfo{author}{A.~Cl{é}ment},
  \bibinfo{author}{G.~Perrin}, \bibinfo{author}{J.~Giorla},
  \bibinfo{author}{H.~Wynn},
\newblock \bibinfo{title}{Group kernels for gaussian process metamodels with
  categorical inputs},
\newblock \bibinfo{journal}{SIAM Journal on Uncertainty Quantification}
  \bibinfo{volume}{8} (\bibinfo{year}{2020}) \bibinfo{pages}{775--806}.
\bibitem[{Garrido-Merchán and Hernández-Lobato(2020)}]{GMHL}
\bibinfo{author}{E.~C. Garrido-Merchán},
  \bibinfo{author}{D.~Hernández-Lobato},
\newblock \bibinfo{title}{Dealing with categorical and integer-valued variables
  in bayesian optimization with gaussian processes},
\newblock \bibinfo{journal}{Neurocomputing} \bibinfo{volume}{380}
  (\bibinfo{year}{2020}) \bibinfo{pages}{20--35}.
\bibitem[{Halstrup(2016)}]{Gower}
\bibinfo{author}{M.~Halstrup}, \bibinfo{title}{Black-Box Optimization of Mixed
  Discrete-Continuous Optimization Problems}, Ph.D. thesis, TU Dortmund,
  \bibinfo{year}{2016}.
\bibitem[{Cuesta-Ramirez et~al.(2021)Cuesta-Ramirez, {Le Riche}, Roustant,
  Perrin, Durantin, and Gliere}]{cuesta2021comparison}
\bibinfo{author}{J.~Cuesta-Ramirez}, \bibinfo{author}{R.~{Le Riche}},
  \bibinfo{author}{O.~Roustant}, \bibinfo{author}{G.~Perrin},
  \bibinfo{author}{C.~Durantin}, \bibinfo{author}{A.~Gliere},
\newblock \bibinfo{title}{A comparison of mixed-variables bayesian optimization
  approaches},
\newblock \bibinfo{journal}{Advanced Modeling and Simulation in Engineering
  Sciences} \bibinfo{volume}{9} (\bibinfo{year}{2021}) \bibinfo{pages}{1--29}.
\bibitem[{Hutter et~al.(2011)Hutter, Hoos, and Leyton-Brown}]{SMAC}
\bibinfo{author}{F.~Hutter}, \bibinfo{author}{H.~Hoos},
  \bibinfo{author}{K.~Leyton-Brown},
\newblock \bibinfo{title}{Sequential model-based optimization for general
  algorithm configuration},
\newblock in: \bibinfo{booktitle}{International Conference on Learning and
  Intelligent Optimization}, \bibinfo{year}{2011}.
\bibitem[{Bergstra et~al.(2011)Bergstra, Bardenet, Bengio, and K\'{e}gl}]{TPE}
\bibinfo{author}{J.~S. Bergstra}, \bibinfo{author}{R.~Bardenet},
  \bibinfo{author}{Y.~Bengio}, \bibinfo{author}{B.~K\'{e}gl},
\newblock \bibinfo{title}{Algorithms for hyper-parameter optimization},
\newblock in: \bibinfo{booktitle}{25th Annual Conference on Neural Information
  Processing Systems}, \bibinfo{year}{2011}.
\bibitem[{Bliek et~al.(2021)Bliek, Guijt, Verwer, and de~Weerdt}]{Relu-surr}
\bibinfo{author}{L.~Bliek}, \bibinfo{author}{A.~Guijt},
  \bibinfo{author}{S.~Verwer}, \bibinfo{author}{M.~de~Weerdt},
\newblock \bibinfo{title}{Black-box mixed-variable optimisation using a
  surrogate model that satisfies integer constraints},
\newblock in: \bibinfo{booktitle}{Proceedings of the Genetic and Evolutionary
  Computation Conference Companion}, \bibinfo{year}{2021}.
\bibitem[{Papalexopoulos et~al.(2022)Papalexopoulos, Tjandraatmadja, Anderson,
  Vielma, and Belanger}]{nn-surr}
\bibinfo{author}{T.~Papalexopoulos}, \bibinfo{author}{C.~Tjandraatmadja},
  \bibinfo{author}{R.~Anderson}, \bibinfo{author}{J.~P. Vielma},
  \bibinfo{author}{D.~Belanger},
\newblock \bibinfo{title}{Constrained discrete black-box optimization using
  mixed-integer programming},
\newblock in: \bibinfo{booktitle}{International Conference on Machine
  Learning}, \bibinfo{year}{2022}.
\bibitem[{Nie and Racine(2012)}]{splines}
\bibinfo{author}{Z.~Nie}, \bibinfo{author}{J.~Racine},
\newblock \bibinfo{title}{The crs package: nonparametric regression splines for
  continuous and categorical predictors},
\newblock \bibinfo{journal}{R Journal} \bibinfo{volume}{4}
  (\bibinfo{year}{2012}) \bibinfo{pages}{48--56}.
\bibitem[{Herrera et~al.(2014)Herrera, Guglielmetti, Xiao, and
  Coelho}]{herrera}
\bibinfo{author}{M.~Herrera}, \bibinfo{author}{A.~Guglielmetti},
  \bibinfo{author}{M.~Xiao}, \bibinfo{author}{R.~F. Coelho},
\newblock \bibinfo{title}{Metamodel-assisted optimization based on multiple
  kernel regression for mixed variables},
\newblock \bibinfo{journal}{Structural and multidisciplinary optimization}
  \bibinfo{volume}{49} (\bibinfo{year}{2014}) \bibinfo{pages}{979--991}.
\bibitem[{Moraglio and Kattan(2011)}]{RBF_geo}
\bibinfo{author}{A.~Moraglio}, \bibinfo{author}{A.~Kattan},
\newblock \bibinfo{title}{Geometric generalisation of surrogate model based
  optimisation to combinatorial spaces},
\newblock in: \bibinfo{booktitle}{Evolutionary Computation in Combinatorial
  Optimization}, \bibinfo{year}{2011}.
\bibitem[{{Munoz Zuniga} and Sinoquet(2020)}]{CAT-EGO}
\bibinfo{author}{M.~{Munoz Zuniga}}, \bibinfo{author}{D.~Sinoquet},
\newblock \bibinfo{title}{Global optimization for mixed categorical-continuous
  variables based on gaussian process models with a randomized categorical
  space exploration step},
\newblock \bibinfo{journal}{INFOR: Information Systems and Operational
  Research} \bibinfo{volume}{58} (\bibinfo{year}{2020})
  \bibinfo{pages}{310--341}.
\bibitem[{Nguyen et~al.(2020)Nguyen, Gupta, Rana, Shilto, and
  Venkatesh}]{Bandit-BO}
\bibinfo{author}{D.~Nguyen}, \bibinfo{author}{S.~Gupta},
  \bibinfo{author}{S.~Rana}, \bibinfo{author}{A.~Shilto},
  \bibinfo{author}{S.~Venkatesh},
\newblock \bibinfo{title}{Bayesian optimization for categorical and
  category-specific continuous inputs},
\newblock in: \bibinfo{booktitle}{AAAI-20 Technical Tracks},
  \bibinfo{year}{2020}.
\bibitem[{Roy et~al.(2019)Roy, Crossley, Stanford, Moore, and Gray}]{AMIEGO}
\bibinfo{author}{S.~Roy}, \bibinfo{author}{W.~A. Crossley},
  \bibinfo{author}{B.~K. Stanford}, \bibinfo{author}{K.~T. Moore},
  \bibinfo{author}{J.~S. Gray},
\newblock \bibinfo{title}{A mixed integer efficient global optimization
  algorithm with multiple infill strategy - applied to a wing topology
  optimization problem},
\newblock in: \bibinfo{booktitle}{AIAA SciTech 2019 Forum},
  \bibinfo{year}{2019}.
\bibitem[{Abramson et~al.(2007)Abramson, Audet, and Dennis}]{Mixed_Abramson}
\bibinfo{author}{M.~A. Abramson}, \bibinfo{author}{C.~Audet},
  \bibinfo{author}{J.~Dennis},
\newblock \bibinfo{title}{Filter pattern search algorithms for mixed variable
  constrained optimization problems},
\newblock \bibinfo{journal}{Pacific Journal of Optimization}
  \bibinfo{volume}{3} (\bibinfo{year}{2007}) \bibinfo{pages}{477–500}.
\bibitem[{Bouhlel et~al.(2019)Bouhlel, Hwang, Bartoli, Lafage, Morlier, and
  Martins}]{SMT2019}
\bibinfo{author}{M.~A. Bouhlel}, \bibinfo{author}{J.~T. Hwang},
  \bibinfo{author}{N.~Bartoli}, \bibinfo{author}{R.~Lafage},
  \bibinfo{author}{J.~Morlier}, \bibinfo{author}{J.~R.~A. Martins},
\newblock \bibinfo{title}{A python surrogate modeling framework with
  derivatives},
\newblock \bibinfo{journal}{Advances in Engineering Software}
  \bibinfo{volume}{135} (\bibinfo{year}{2019}) \bibinfo{pages}{102662}.
\bibitem[{Kr\"{u}gener et~al.(2022)Kr\"{u}gener, {Zapata Usandivaras}, ,
  Bauerheim, and Urbano}]{DL1}
\bibinfo{author}{M.~Kr\"{u}gener}, \bibinfo{author}{J.~{Zapata Usandivaras}}, ,
  \bibinfo{author}{M.~Bauerheim}, \bibinfo{author}{A.~Urbano},
\newblock \bibinfo{title}{Coaxial-injector surrogate modeling based on
  reynolds-averaged navier–stokes simulations using deep learning},
\newblock \bibinfo{journal}{Journal of Propulsion and Power}
  \bibinfo{volume}{38} (\bibinfo{year}{2022}) \bibinfo{pages}{783--798}.
\bibitem[{Li et~al.(2022)Li, Zhang, Tay, Liu, Cui, Chew, and Khoo}]{DL2}
\bibinfo{author}{J.~Li}, \bibinfo{author}{M.~Zhang}, \bibinfo{author}{C.~M.~J.
  Tay}, \bibinfo{author}{N.~Liu}, \bibinfo{author}{Y.~Cui},
  \bibinfo{author}{S.~C. Chew}, \bibinfo{author}{B.~C. Khoo},
\newblock \bibinfo{title}{Low-reynolds-number airfoil design optimization using
  deep-learning-based tailored airfoil modes},
\newblock \bibinfo{journal}{Aerospace Science and Technology}
  \bibinfo{volume}{121} (\bibinfo{year}{2022}) \bibinfo{pages}{107309}.
\bibitem[{Li et~al.(2020)Li, Zhang, Martins, and Shu}]{DL3}
\bibinfo{author}{J.~Li}, \bibinfo{author}{M.~Zhang}, \bibinfo{author}{J.~R.
  R.~A. Martins}, \bibinfo{author}{C.~Shu},
\newblock \bibinfo{title}{Efficient aerodynamic shape optimization with
  deep-learning-based geometric filtering},
\newblock \bibinfo{journal}{AIAA Journal} \bibinfo{volume}{58}
  (\bibinfo{year}{2020}) \bibinfo{pages}{4243--4259}.
\bibitem[{{Zapata Usandivaras} et~al.(2022){Zapata Usandivaras}, Bauerheim,
  Benedicte, and Urbano}]{DL4}
\bibinfo{author}{J.~{Zapata Usandivaras}}, \bibinfo{author}{M.~Bauerheim},
  \bibinfo{author}{C.~Benedicte}, \bibinfo{author}{A.~Urbano},
\newblock \bibinfo{title}{Large eddy simulations and deep learning for the
  investigation of recess variation of a shear-coaxial injector},
\newblock in: \bibinfo{booktitle}{Space Propulsion Conference 2022},
  \bibinfo{year}{2022}.
\bibitem[{Ming et~al.(2022)Ming, Williamson, and Guillas}]{DGP1}
\bibinfo{author}{D.~Ming}, \bibinfo{author}{D.~Williamson},
  \bibinfo{author}{S.~Guillas},
\newblock \bibinfo{title}{Deep gaussian process emulation using stochastic
  imputation},
\newblock \bibinfo{journal}{Technometrics} \bibinfo{volume}{0}
  (\bibinfo{year}{2022}) \bibinfo{pages}{1--12}.
\bibitem[{Izzaturrahman et~al.(2021)Izzaturrahman, Palar, Zuhal, and
  Shimoyama}]{DGP2}
\bibinfo{author}{M.~F. Izzaturrahman}, \bibinfo{author}{P.~S. Palar},
  \bibinfo{author}{L.~Zuhal}, \bibinfo{author}{K.~Shimoyama},
\newblock \bibinfo{title}{Modeling non-stationarity with deep gaussian
  processes: Applications in aerospace engineering},
\newblock in: \bibinfo{booktitle}{AIAA SciTech 2022 Forum},
  \bibinfo{year}{2021}.
\bibitem[{Forrester et~al.(2008)Forrester, Sobester, and Keane}]{forrester}
\bibinfo{author}{A.~Forrester}, \bibinfo{author}{A.~Sobester},
  \bibinfo{author}{A.~Keane}, \bibinfo{title}{Engineering Design via Surrogate
  Modelling: A Practical Guide}, \bibinfo{publisher}{Wiley},
  \bibinfo{year}{2008}.
\bibitem[{Duvenaud(2014)}]{GP14}
\bibinfo{author}{D.~Duvenaud}, \bibinfo{title}{Automatic model construction
  with Gaussian processes}, Ph.D. thesis, University of Cambridge,
  \bibinfo{year}{2014}.
\bibitem[{Rossi(2018)}]{MLE}
\bibinfo{author}{R.~J. Rossi}, \bibinfo{title}{Mathematical statistics: an
  introduction to likelihood based inference}, \bibinfo{publisher}{John Wiley
  \& Sons}, \bibinfo{year}{2018}.
\bibitem[{{De L{\'a}zaro} et~al.(2015){De L{\'a}zaro}, Moreno, Santiago, and
  {da Silva Neto}}]{de2015optimizing}
\bibinfo{author}{J.~M.~B. {De L{\'a}zaro}}, \bibinfo{author}{A.~P. Moreno},
  \bibinfo{author}{O.~L. Santiago}, \bibinfo{author}{A.~J. {da Silva Neto}},
\newblock \bibinfo{title}{Optimizing kernel methods to reduce dimensionality in
  fault diagnosis of industrial systems},
\newblock \bibinfo{journal}{Computers \& Industrial Engineering}
  \bibinfo{volume}{87} (\bibinfo{year}{2015}) \bibinfo{pages}{140--149}.
\bibitem[{{De L{\'a}zaro} et~al.(2022){De L{\'a}zaro}, {Cruz-Corona}, , {da
  Silva Neto}, and Santiago}]{bernal2022criteria}
\bibinfo{author}{J.~M.~B. {De L{\'a}zaro}}, \bibinfo{author}{C.~{Cruz-Corona}},
  , \bibinfo{author}{A.~J. {da Silva Neto}}, \bibinfo{author}{O.~L. Santiago},
\newblock \bibinfo{title}{Criteria for optimizing kernel methods in fault
  monitoring process: A survey},
\newblock \bibinfo{journal}{ISA transactions} \bibinfo{volume}{127}
  (\bibinfo{year}{2022}) \bibinfo{pages}{259--272}.
\bibitem[{Lee(2011)}]{Lee2011}
\bibinfo{author}{H.~Lee}, \bibinfo{title}{Gaussian Processes},
  volume~\bibinfo{volume}{5}, \bibinfo{publisher}{Springer Berlin Heidelberg},
  \bibinfo{year}{2011}, pp. \bibinfo{pages}{575--577}.
\bibitem[{Golovin et~al.(2017)Golovin, Solnik, Moitra, Kochanski, Karro, and
  Sculley}]{one-hot}
\bibinfo{author}{D.~Golovin}, \bibinfo{author}{B.~Solnik},
  \bibinfo{author}{S.~Moitra}, \bibinfo{author}{G.~Kochanski},
  \bibinfo{author}{J.~Karro}, \bibinfo{author}{D.~Sculley},
\newblock \bibinfo{title}{Google vizier: A service for black-box optimization},
\newblock in: \bibinfo{booktitle}{Proceedings of the 23rd ACM SIGKDD
  International Conference on Knowledge Discovery and Data Mining},
  \bibinfo{year}{2017}.
\bibitem[{Rebonato and Jaeckel(2001)}]{HS}
\bibinfo{author}{R.~Rebonato}, \bibinfo{author}{P.~Jaeckel},
\newblock \bibinfo{title}{The most general methodology to create a valid
  correlation matrix for risk management and option pricing purposes},
\newblock \bibinfo{journal}{Journal of Risk} \bibinfo{volume}{2}
  (\bibinfo{year}{2001}) \bibinfo{pages}{17–27}.
\bibitem[{Rapisarda et~al.(2007)Rapisarda, Brigo, and Mercurio}]{HS_Jacobi}
\bibinfo{author}{F.~Rapisarda}, \bibinfo{author}{D.~Brigo},
  \bibinfo{author}{F.~Mercurio},
\newblock \bibinfo{title}{Parameterizing correlations: a geometric
  interpretation},
\newblock \bibinfo{journal}{IMA Journal of Management Mathematics}
  \bibinfo{volume}{18} (\bibinfo{year}{2007}) \bibinfo{pages}{55--73}.
\bibitem[{Qian et~al.(2008)Qian, Wu, and Wu}]{PDUDE}
\bibinfo{author}{P.~Z.~G. Qian}, \bibinfo{author}{H.~Wu},
  \bibinfo{author}{C.~F.~J. Wu},
\newblock \bibinfo{title}{Gaussian process models for computer experiments with
  qualitative and quantitative factors},
\newblock \bibinfo{journal}{Technometrics} \bibinfo{volume}{50}
  (\bibinfo{year}{2008}) \bibinfo{pages}{383--396}.
\bibitem[{Hadamard(1910)}]{hadamard}
\bibinfo{author}{J.~Hadamard}, \bibinfo{title}{Sur quelques applications de
  l'indice de Kronecker}, \bibinfo{publisher}{Bussiere}, \bibinfo{year}{1910}.
\bibitem[{Bapat and Raghavan(1997)}]{schurmatapp}
\bibinfo{author}{R.~B. Bapat}, \bibinfo{author}{T.~E.~S. Raghavan},
  \bibinfo{title}{Nonnegative Matrices and Applications},
  \bibinfo{publisher}{Cambridge University Press}, \bibinfo{year}{1997}.
\bibitem[{Schoenberg(1938)}]{Schoenberg}
\bibinfo{author}{I.~J. Schoenberg},
\newblock \bibinfo{title}{Metric spaces and positive definite functions},
\newblock \bibinfo{journal}{Transactions of the American Mathematical Society}
  \bibinfo{volume}{44} (\bibinfo{year}{1938}) \bibinfo{pages}{522--536}.
\bibitem[{Horn and Johnson(2012)}]{horn2012matrix}
\bibinfo{author}{R.~A. Horn}, \bibinfo{author}{C.~R. Johnson},
  \bibinfo{title}{Matrix analysis}, \bibinfo{publisher}{Cambridge university
  press}, \bibinfo{year}{2012}.
\bibitem[{Vilenkin and Klimyk(1995)}]{hypersphere}
\bibinfo{author}{N.~Y. Vilenkin}, \bibinfo{author}{A.~U. Klimyk},
  \bibinfo{title}{Representations of Lie groups and special functions},
  \bibinfo{publisher}{Springer}, \bibinfo{year}{1995}.
\bibitem[{Powell(1994)}]{COBYLA}
\bibinfo{author}{M.~J.~D. Powell}, \bibinfo{title}{A direct search optimization
  method that models the objective and constraint functions by linear
  interpolation}, volume \bibinfo{volume}{275}, \bibinfo{publisher}{Springer},
  \bibinfo{year}{1994}, pp. \bibinfo{pages}{51--67}.
\bibitem[{Jin et~al.(2005)Jin, Chen, and Sudjianto}]{LHS}
\bibinfo{author}{R.~Jin}, \bibinfo{author}{W.~Chen},
  \bibinfo{author}{A.~Sudjianto},
\newblock \bibinfo{title}{An efficient algorithm for constructing optimal
  design of computer experiments},
\newblock \bibinfo{journal}{Journal of Statistical Planning and Inference}
  \bibinfo{volume}{2} (\bibinfo{year}{2005}) \bibinfo{pages}{545--554}.
\bibitem[{Demay et~al.(2022)Demay, Iooss, Gratiet, and Marrel}]{PVA}
\bibinfo{author}{C.~Demay}, \bibinfo{author}{B.~Iooss}, \bibinfo{author}{L.~L.
  Gratiet}, \bibinfo{author}{A.~Marrel},
\newblock \bibinfo{title}{Model selection based on validation criteria for
  gaussian process regression: An application with highlights on the predictive
  variance},
\newblock \bibinfo{journal}{Quality and Reliability Engineering International}
  \bibinfo{volume}{38} (\bibinfo{year}{2022}) \bibinfo{pages}{1482--1500}.
\bibitem[{Cheng et~al.(2015)Cheng, Younis, Hajikolaei, and
  Wang}]{Cheng2015TrustRB}
\bibinfo{author}{G.~H. Cheng}, \bibinfo{author}{A.~Younis},
  \bibinfo{author}{K.~H. Hajikolaei}, \bibinfo{author}{G.~G. Wang},
\newblock \bibinfo{title}{Trust region based mode pursuing sampling method for
  global optimization of high dimensional design problems},
\newblock \bibinfo{journal}{Journal of Mechanical Design} \bibinfo{volume}{137}
  (\bibinfo{year}{2015}) \bibinfo{pages}{021407}.
\bibitem[{Oune and Bostanabad(2021)}]{oune2021latent}
\bibinfo{author}{N.~Oune}, \bibinfo{author}{R.~Bostanabad},
\newblock \bibinfo{title}{Latent map gaussian processes for mixed variable
  metamodeling},
\newblock \bibinfo{journal}{Computer Methods in Applied Mechanics and
  Engineering} \bibinfo{volume}{387} (\bibinfo{year}{2021})
  \bibinfo{pages}{114128}.
\bibitem[{Schmollgruber et~al.(2019)Schmollgruber, D{\"o}ll, Hermetz, Liaboeuf,
  Ridel, Cafarelli, Atinault, Fran{\c c}ois, and Paluch}]{schmollgruber}
\bibinfo{author}{P.~Schmollgruber}, \bibinfo{author}{C.~D{\"o}ll},
  \bibinfo{author}{J.~Hermetz}, \bibinfo{author}{R.~Liaboeuf},
  \bibinfo{author}{M.~Ridel}, \bibinfo{author}{I.~Cafarelli},
  \bibinfo{author}{O.~Atinault}, \bibinfo{author}{C.~Fran{\c c}ois},
  \bibinfo{author}{B.~Paluch},
\newblock \bibinfo{title}{Multidisciplinary exploration of {DRAGON}: an {ONERA}
  hybrid electric distributed propulsion concept},
\newblock in: \bibinfo{booktitle}{AIAA SciTech 2019 Forum},
  \bibinfo{year}{2019}.
\bibitem[{David et~al.(2021)David, Delbecq, Defoort, Schmollgruber, Benard, and
  Pommier-Budinger}]{David_2021}
\bibinfo{author}{C.~David}, \bibinfo{author}{S.~Delbecq},
  \bibinfo{author}{S.~Defoort}, \bibinfo{author}{P.~Schmollgruber},
  \bibinfo{author}{E.~Benard}, \bibinfo{author}{V.~Pommier-Budinger},
\newblock \bibinfo{title}{From {FAST} to {FAST}-{OAD}: An open source framework
  for rapid overall aircraft design},
\newblock \bibinfo{journal}{{IOP} Conference Series: Materials Science and
  Engineering} \bibinfo{volume}{1024} (\bibinfo{year}{2021})
  \bibinfo{pages}{012062}.
\bibitem[{Bouhlel et~al.(2018)Bouhlel, Bartoli, Regis, Otsmane, and
  Morlier}]{Bouhlel18}
\bibinfo{author}{M.~A. Bouhlel}, \bibinfo{author}{N.~Bartoli},
  \bibinfo{author}{R.~Regis}, \bibinfo{author}{A.~Otsmane},
  \bibinfo{author}{J.~Morlier},
\newblock \bibinfo{title}{Efficient global optimization for high-dimensional
  constrained problems by using the kriging models combined with the partial
  least squares method},
\newblock \bibinfo{journal}{Engineering Optimization} \bibinfo{volume}{50}
  (\bibinfo{year}{2018}) \bibinfo{pages}{2038--2053}.

\end{thebibliography}

\end{document}